\documentclass{amsart}

\usepackage{amsfonts,amssymb,amscd,amsmath,enumerate,xypic}
\xyoption{curve}

\newcommand{\qnd}{\quad\text{and}\quad}

\newcommand{\foral}{\qquad\text{for all}\quad}

\newcommand{\where}{\qquad\text{where}\quad}
\newcommand{\with}{\qquad\text{with}\quad}

\newcommand{\ssstyle}{\scriptscriptstyle}

\newcommand{\les}{\operatorname{\ssstyle\leqslant}}

\newcommand{\sshift}{{\operatorname{\scriptstyle\Sigma}}}
\newcommand{\shift}[2]{\sshift^{#1}#2}

\newcommand{\sign}[2]{(-1)^{|#1||#2|}}

\newcommand{\col}{\colon}
\newcommand{\dd}{\partial}

\newcommand{\wt}{\widetilde}
\newcommand{\wh}{\widehat}

\newcommand{\arto}[1]{\overset{#1}{\longrightarrow}}

\newcommand{\tra}{\twoheadrightarrow}

\newcommand{\quism}{\text{quasi\-isomorphism}}

\newcommand{\perf}{{small}}
\newcommand{\vsmall}{{virtually\- small}}
\newcommand{\psmall}{{proxy\- small}}

\newcommand{\card}{\operatorname{card}}
\newcommand{\id}{\operatorname{id}}

\newcommand{\Ker}{\operatorname{Ker}}

\newcommand{\codim}{\operatorname{codim}}

\newcommand{\hilb}[2]{\operatorname{Hilb}_{#1}(#2)}

\newcommand{\rank}{\operatorname{rank}}

\newcommand{\res}[1]{{\boldsymbol k}(#1)}
\newcommand{\spec}{\operatorname{Spec}}
\newcommand{\supp}{\operatorname{supp}}
\newcommand{\Supp}{\operatorname{Supp}}
\newcommand{\V}{{\operatorname{V}}}

\newcommand{\ZZ}{{\operatorname{Z}}}

\newcommand{\dtensor}[1]{{\otimes^{\mathbf L}_{#1}}}
\newcommand{\ext}[4]{\operatorname{Ext}^{#1}_{#2}\left(#3,#4\right){}}
\newcommand{\Hom}[3]{\operatorname{Hom}_{#1}\left(#2,#3\right){}}
\newcommand{\hh}[1]{\operatorname{H}(#1)}
\newcommand{\HH}[2]{\operatorname{H}_{#1}(#2)}
\newcommand{\rhom}[3]{\operatorname{{\mathbf R}Hom}_{#1}\left(#2,#3\right){}}
\newcommand{\tor}[4]{\operatorname{Tor}_{#1}^{#2}\left(#3,#4\right){}}

\newcommand{\cxy}[2]{\operatorname{cx}_{#1}{#2}}
\newcommand{\curv}[2]{\operatorname{curv}_{#1}{#2}}
\newcommand{\injcxy}[2]{\operatorname{injcx}_{#1}{#2}}
\newcommand{\injcurv}[2]{\operatorname{injcurv}_{#1}{#2}}
\newcommand{\fdim}{\operatorname{fd}}
\newcommand{\gdim}{\operatorname{G-dim}}
\newcommand{\idim}{\operatorname{id}}
\newcommand{\pdim}{\operatorname{pd}}

\newcommand{\amp}{\operatorname{amp}}
\newcommand{\builds}{{\Longrightarrow}}
\newcommand{\rbuilds}[1]{\underset{#1}{\Longrightarrow}}

\newcommand{\rthick}[2]{{\mathsf{Thick}}_{#1}[#2]}

\newcommand{\comp}[2]{{\mathcal C}^{#1}\!{(#2)}}

\newcommand{\diag}[2]{{#2}^{e}}
\newcommand{\koszul}[2]{\operatorname{K}[#1\,;#2]}

\newcommand{\lcat}[1]{\mathcal{D}\left(#1\right)}
\newcommand{\exact}[3]{{#1\longrightarrow #2 \longrightarrow #3 \longrightarrow \sshift #1}}
\newcommand{\Lch}[2]{\mathbf{R}\Gamma_{#1}\left(#2\right)}
\newcommand{\Lho}[2]{\mathbf{L}\Lambda^{#1}\left(#2\right)}
\newcommand{\subcat}[2]{{\mathcal{D}^{#1}\!\left(#2\right)}}
\newcommand{\colim}[1]{{\scriptstyle\coprod}_{#1}}

\newcommand{\bass}[1]{\operatorname{I}^{#1}(t)}
\newcommand{\rbass}[2]{\operatorname{I}_{#1}^{#2}(t)}
\newcommand{\poin}[1]{\operatorname{P}_{#1}(t)}
\newcommand{\rpoin}[2]{\operatorname{P}^{#1}_{#2}(t)}

\newcommand{\bsq}{{\boldsymbol q}}

\newcommand{\bsx}{{\boldsymbol x}}
\newcommand{\bsy}{{\boldsymbol y}}

\newcommand{\BN}{{\mathbb N}}

\newcommand{\BR}{{\mathbb R}}
\newcommand{\BZ}{{\mathbb Z}}

\newcommand{\calf}{{\mathcal F}}
\newcommand{\calt}{{\mathcal T}}

\newcommand{\fm}{{\mathfrak m}}

\newcommand{\fn}{{\mathfrak n}}
\newcommand{\fp}{{\mathfrak p}}
\newcommand{\fq}{{\mathfrak q}}

\newcommand{\eps}{{\epsilon}}

\newcommand{\vf}{{\varphi}}

\swapnumbers

\theoremstyle{plain}
\newtheorem{theorem}{Theorem}[section]
\newtheorem*{Theorem}{Theorem}

\newtheorem{itheorem}{Theorem}

\newtheorem{icorollary}[itheorem]{Corollary}

\newtheorem{proposition}[theorem]{Proposition}

\newtheorem{lemma}[theorem]{Lemma}

\newtheorem{corollary}[theorem]{Corollary}

\newtheorem{itchunk}[theorem]{}

\theoremstyle{definition}

\newtheorem{example}[theorem]{Example}

\newtheorem{problem}[theorem]{Problem}

\newenvironment{bfchunk}{\begin{chunk}\textbf}{\end{chunk}}
\newtheorem{chunk}[theorem]{}

\theoremstyle{remark}

\newtheorem*{Claim}{Claim}

\newtheorem{remark}[theorem]{Remark}
\newtheorem{remarks}[theorem]{Remarks}

\newtheorem{question}[theorem]{Question}
\newtheorem*{Question}{Question}

\begin{document}

\title[Finiteness in derived categories] {Finiteness in derived categories of local rings}

\author{W.~Dwyer} 
\address{University of Notre Dame, Notre Dame, IN
  46556, U.S.A.}
\email{dwyer.1@nd.edu} 

\author{J.~P.~C.~Greenlees} 
 \address{Pure Mathematics,
  Hicks Building, University of Sheffield, S3 7RH, U.K.}
\email{j.greenlees@shef.ac.uk}

\author{S.~Iyengar} \address{Department of Mathematics,
  University of Nebraska, Lincoln, NE 68588, U.S.A.}
 \email{iyengar@math.unl.edu}



\thanks {S.I. was partly supported by the NSF under grant DMS 0302892.}

\begin{abstract}
  New homotopy invariant finiteness conditions on modules over commutative rings are
  introduced, and their properties are studied systematically. A number of finiteness
  results for classical homological invariants like flat dimension, injective dimension,
  and Gorenstein dimension, are established. It is proved that these specialize to give
  results concerning modules over complete intersection local rings. A noteworthy feature
  is the use of techniques based on thick subcategories of derived categories.
\end{abstract}

\maketitle


\tableofcontents  

\setcounter{section}{0}

\section{Introduction}
This paper investigates several homotopy invariant finiteness conditions on modules and
complexes over commutative noetherian rings. This leads to new results in classical
commutative algebra, some of which are discussed in the introduction.  First, we describe
our philosophy.

Our earlier work \cite{DGI} transplants ideas from commutative algebra into a homotopy
theoretic context. The process is conceptually illuminating since it shows that several
classical duality theorems in algebra and topology are not just analogous but actually
manifestations of the same phenomenon.  However, the main justification for the approach
is that it leads to interesting and concrete new results in homotopy theory. In the
present paper we reverse the process, and use ideas inspired by homotopy theory to solve
problems in ring theory.

Before one can apply a concept from algebra in the topological context one has to express
it in a homotopy invariant form. The idea of formulating statements about a ring in terms
of its derived category is a familiar one, but because the framework of \cite{DGI} was
more general than that, we were forced to concentrate on the crudest and most robust
properties. It is striking that these give new insights even in commutative ring theory.
The homological characterization of regular local rings illustrates this point.

A local ring $R$ is said to be regular if its maximal ideal $\fm$ is generated by a
regular sequence.  Auslander, Buchsbaum, and Serre proved that $R$ is regular if and only
if the residue field $k$ has finite projective dimension.  Since the finiteness of the
projective dimension of a homologically finite complex $M$ is detected by $\ext{*}RMk$,
regularity is detected at the level of derived categories.  However, one can reformulate
it in structural terms: $k$ has finite projective dimension precisely when it is
equivalent in the derived category to a perfect complex, that is to say, a finite complex
consisting of finitely generated projectives. Furthermore, perfectness can be captured
entirely categorically: $M$ is perfect if and only if it is \emph{small}, in the sense
that $\rhom RM-$ preserves arbitrary direct sums in the derived category.  Thus one
reaches the homotopy invariant formulation: $R$ is regular if and only if $k$ is small in
the derived category. The selection of the particular module $k$ is unnecessary: \emph{$R$
  is regular if and only if every homologically finite complex $M$ is small in the derived
  category}.  The advantages of the homological formulation are well known: for example,
it allows one to prove that localizations of regular local rings are regular. The
structural formulation shows that regularity is preserved by equivalences of derived
categories, and makes sense even when one has no notion of elements or of Ext.

Small complexes remain important in commutative ring theory, but the focus of our work is
on \emph{\vsmall} complexes: those complexes $M$ with the much weaker property that either
$\hh M=0$, or the thick subcategory generated by $M$ contains a non-trivial small
complex. A detailed explanation of the terms involved in the preceding definition,
including some background on thick subcategories, is provided in Section
\ref{section:Thicksubcategories}.  We study also the closely related \emph{\psmall}
complexes, which played a crucial role in \cite{DGI} and provided the original motivation
for the present investigations.

Evidently every \perf\, complex is \vsmall.  However there are many more virtually small
complexes than small ones; for instance, if $M$ is \vsmall, then so is $M\oplus X$, for
any complex $X$ of $R$-modules, as long as $\hh M\ne 0$.  A better illustration is the
result below, which contains parts of Theorems \eqref{vsmall:almostfinite} and
\eqref{vsmall:regulardescent}. Given a complex $N$ of $S$-modules, we say $\hh N$ is
finite (respectively, artinian) if the $S$-module $\HH iN$ is finite (respectively,
artinian) for each $i$, and zero for $|i|\gg 0$.

\begin{itheorem}
  Let $\vf\col R \to S$ be a local homomorphism and $N$ a complex  of
  $S$-modules with $\hh N$ either finite or artinian.  If either
  $\fdim_RN$ is finite  or $S$ is regular, then $N$ is \vsmall\, over $R$.
\end{itheorem}

Specializing $\vf$ to the surjective homomorphism $R\to k$ and $N$ to $k$, we see that the
residue field $k$ is \vsmall\, over $R$; alternatively, it is not hard to check that the
Koszul complex of $R$ is in the thick category generated by $k$. This may be a little
unsettling: from the classical point of view, exemplified by the Auslander-Buchsbaum-Serre
theorem, the homological properties of $k$ differ drastically from those of \perf\,
complexes.  Nonetheless, \vsmall\, complexes are very useful as test objects for
finiteness of homological invariants. Sections \ref{section:testobjects} and
\ref{section:Gdimensions} establish a number of results to validate this claim. We now
provide some illustrations.

Let $\psi\col Q\to R$ be a local homomorphism and $M$ a complex of $R$-modules, which is
viewed as a complex of $Q$-modules via $\psi$. It is classical that when $M$ is \perf\,
many homological properties of $R$ over $Q$ are inherited by $M$; for example, it follows
from a result of Cartan-Eilenberg~\cite{CE} that, under this condition on $M$, if
$\fdim_QR$ is finite, then so is $\fdim_QM$.  Recent research, notably that of
Apassov~\cite{Ap}, Foxby and Iyengar~\cite{FI}, and Iyengar and Sather-Wagstaff~\cite{IW},
has uncovered certain situations when the converse is true.  The main discovery of this
work is these are instances of the same phenomenon, which is valid in much greater
generality, and even when $M$ is only assumed to be \vsmall.

For example we may consider flat dimension.  As a converse to the result of
Cartan-Eilenberg recalled above, Foxby and Iyengar~\cite{FI} prove that when $M$ is \perf,
$\fdim_QM$ finite implies $\fdim_QR$ finite.  Furthermore, if $\fdim_Qk$ is finite, then
$Q$ is regular, and hence $\fdim_QR$ is finite; here one is using the
Auslander-Buchsbaum-Serre theorem. We generalize these results as follows:

\begin{itheorem}
  Let $M$ be a \vsmall\, complex of $R$-modules with $\hh M\ne 0$.  If $\fdim_QM$ is
  finite, then $\fdim_QR$ is finite as well.
\end{itheorem}

The proof in \cite{FI} proof is involved, and based on properties of depth for complexes.
In contrast, the one given here is straightforward; it is based on work of
Hopkins~\cite{Ho} and Neeman~\cite{Ne} on thick subcategories of derived categories of
commutative noetherian rings. It was this argument, which is crystallized in
\eqref{proofprinciple}, that convinced us of the relevance of \vsmall\, complexes to
classical commutative ring theory.

We prove numerous other results that demonstrate the efficacy of this technique: Theorems
\eqref{gor:char} and \eqref{test:gdim} have the same flavour as the one above, and deal
with injective dimension and Gorenstein dimension, respectively.  Theorem
\eqref{test:asymptotes} concerns asymptotic invariants, and says that the complexity and
curvature of $R$ (over $\psi$) are smaller than those for any virtually small complex. In
a different direction, Theorem \eqref{amplitude} states that if for a complex $Y$ of
$R$-modules with $\hh Y$ degreewise finite, the amplitude of $M \dtensor RY$ is finite for a
\vsmall\, complex $M$ with $\hh M\ne 0$, then the amplitude of $Y$ is finite.

With the importance of virtual smallness established, Sections \ref{section:Descent} and
\ref{section:Ascent} describe other sources of \vsmall\, complexes and investigate their
behaviour under standard constructions: virtual smallness is much more robust than \perf
ness.

Not everything homologically finite is \vsmall . One example is described in
\eqref{proof:nonexample}, but it is difficult to decide when a complex, even a
homologically finite one, is \vsmall. This is in contrast with the case of \perf\,
complexes where methods are implemented in computer algebra packages such as
\emph{Macaulay}.  It is therefore natural to ask: \emph{over which local rings are all
  homologically finite complexes of modules \vsmall ?}

Our principal result in this direction is Theorem \eqref{vsmall:ciring}, part of which reads:

\begin{itheorem}
  If a local ring $R$ is a complete intersection, then any complex of $R$-modules which is
  homologically either finite or artinian is \vsmall.
\end{itheorem}

Section \ref{section:lci} also contains the necessary background on complete intersection
rings. In view of this theorem, all the results in Sections \ref{section:testobjects} and
\ref{section:Gdimensions} specialize to the case where $R$ is a complete intersection. We
mention one, which is part of Theorem \eqref{ci:descent}. It extends, and its proof
provides a different perspective on, \cite[(R, CI)]{Ap}, which treats the special case
where $M$ is a finite $R$-module. Note that, with $M=R$ one recovers the well-known result
that both regularity and the complete intersection property descend along homomorphisms of
finite flat dimension.

\begin{icorollary}
  Let $Q\to R$ be a local homomorphism and $M$ a complex of $R$-modules
  with $\hh M\ne 0$ either finite or artinian.  If $R$ is a complete intersection
  and $\fdim_QM$ is finite, then $Q$ is a complete intersection and $\codim Q\leq \codim
  R$.
\end{icorollary}

Crucial in the proof of Theorem III are stability properties of virtual smallness. We
consider a ring homomorphism $\vf\col R \to S$ and virtual smallness over $R$ and over $S$
for complexes of $R$-modules.  Following the variance on prime spectra, we say that
virtual smallness \emph{descends along $\vf$} if each complex of $S$-modules that is
virtually small over $S$ is also virtually small over $R$. Most of our discussion on
descent is in Section \ref{section:Descent}, and the following excerpt from Theorem
\eqref{psmall:transitivity} gives a flavour of some of the results found there.

\begin{itheorem}
  Let $\vf\col R\to S$ be a surjective homomorphism of noetherian rings.  If $S $ is
  virtually small as an $R$-module then virtual smallness descends along $\vf$.
\end{itheorem}

Another noteworthy situation when virtual smallness descends is the case when $S$ is
regular and module finite over $R$. In Section \ref{section:Ascent}, we turn to ascent:
virtual smallness \emph{ascends along $\psi\col Q\to R$} if any complex of $R$-modules that is
virtually small over $Q$ is virtually small over $R$.  Ascent is a lot more delicate than
descent, and our main result dealing with this issue is as follows:

\begin{itheorem}
  Let $\psi\col Q \to R$ be a homomorphism of rings. If $R\dtensor QR$ can be finitely
  built from $R$ as a bimodule, then virtual smallness ascends along $\psi$. 
\end{itheorem}

This criterion for ascent is part of Theorem \eqref{psmall:diagonal}, see also
\eqref{diagonal:remark}, and it is a crucial ingredient in our proof of Theorem III.  The
statement itself is explained in detail in Section \ref{section:Ascent}. One reason for
quoting it here is to point out that it is reminiscent of the classical
Hochschild-Kostant-Rosenberg criterion for smoothness: a homomorphism $\psi$ essentially
of finite type is smooth if and only if it is flat and the the $(R\otimes_QR)$-module $R$
is \perf.

A partial converse on Theorem III above is established in Theorem \eqref{cichar:rings}.

\begin{itheorem}
  Let $R$ be local ring. If every homologically finite complex of $R$-modules is \vsmall,
  then $R$ is Gorenstein.
\end{itheorem}

This shows that the class of rings over which every homologically finite complex is
\vsmall\, is sandwiched between complete intersections and Gorenstein rings. It is a
strict subset of the class of Gorenstein rings as shown by Examples
\eqref{nonci:idealization} and \eqref{nonci:macaulay}.  If it coincides with the class of
complete intersections, then this would give a homotopical characterization in terms of
virtual smallness, analogous to the Auslander-Buchsbaum-Serre characterization of
regularity in terms of \perf ness.

A counter-example would be interesting too for it would signal the presence of an
interesting---the results in this article testify to this claim---new family of local
rings. This aspect of our work is in the same spirit as that of Jorgensen and {\c S}ega
\cite{JS}, who have discovered various classes of local rings strictly contained between
complete intersections and Gorenstein rings, which nevertheless retain some traces of the
homological properties of the former. These classes are not directly related to the one
considered here; Examples \eqref{nonci:idealization} and \eqref{nonci:macaulay} make that
clear.

A comment on the exposition: the techniques we use are a mixture of homotopy theory and
commutative algebra, and since we hope to interest both algebraists and topologists, we
have made a special effort has to provide relevant background when appropriate.

\section{Complexes}
Let $R$ be a ring; the standing hypothesis in this article is that rings are commutative.  Let
$M$ be a complex of $R$-modules.  Our labeling convention is homological:
\[
M:= \cdots \to M_{i+1} \arto{\dd_{i+1}} M_i\arto{\dd_i} M_{i-1}\to \cdots
\]
We set $\inf M = \inf \{i\mid \HH iM\ne 0\}$ and $\sup M = \sup \{i\mid \HH iM\ne 0\}$; in
particular, if $\hh M=0$, then $\inf M = \infty$ and $\sup M = -\infty$.  When $\inf
M>-\infty$, respectively, $\sup M<\infty$, one says that $\hh M$ is \emph{bounded below},
respectively, \emph{bounded above}.  The \emph{amplitude} of $M$ is the number
\begin{gather*}
  \amp M = \sup M - \inf M\,;
\end{gather*}
when this number is finite, then $\hh M$ is said to be \emph{bounded}.

If for each integer $i$ the $R$-module $\HH iM$ is finite, that is to say, finitely
generated, then $\hh M$ is \emph{degreewise finite}.  When $\hh M$ is both degreewise
finite and bounded, we say that $\hh M$ is \emph{finite}, or that \emph{$M$ is
  homologically finite}. Sometimes, we consider also the property that for each integer
$i$, the $R$-module $\HH iM$ is artinian; in other words, $\hh M$ is \emph{degreewise
  artinian}. When $\hh M$ is also bounded, we drop the qualifier `degreewise', or write:
\emph{$M$ is homologically artinian}.

Next, we provide a short resume of the homological algebra of complexes, referring the
reader to the work of Avramov and Foxby \cite{AF:hd} on homological dimensions for greater
detail. Recall that a morphism $\alpha$ of complexes is a \emph{\quism} if $\hh
\alpha$ is bijective.

\begin{bfchunk}{Resolutions.}
  A complex $P$ of $R$-modules is \emph{projective} when it has the property that for each
morphism $\alpha\col M\to N$ of complexes, if $\alpha$ is surjective and a \quism,
  the so is
\[
\Hom RP\alpha\col \Hom RPM \to \Hom RPN\,;
\]
in other words, $P$ preserves surjective {\quism s}.  In \cite{AF:hd}, such a $P$ is said
to be DG projective, in order to distinguish it from other flavours of projectives; we
omit the `DG' for this is the only kind of projective complexes considered here.

Each complex $M$ of $R$-modules has a \emph{projective resolution}: a {\quism} $P\to M$
with $P$ a projective complex. When $R$ is noetherian and $M$ is homologically finite and
bounded below, $P$ may be so constructed that each $P_i$ is finitely generated and $P_i=0$
for $i<\inf M$.  One has also \emph{flat} complexes and \emph{injective} complexes, and
the corresponding resolutions of $M$.
\end{bfchunk}

\begin{bfchunk}{Derived functors.}
  The derived category of $R$ is denoted $\lcat R$; when in it the notation $\simeq$
  indicates an isomorphism; we recommend Hartshorne \cite{Ha} and Verdier \cite{Ve} as
  basic references on derived categories and derived functors.
  
  Given complexes $M$ and $N$ of $R$-modules, $M\dtensor RN$ and $\rhom RMN$ denote the
  derived functor of the tensor product functor and the homomorphisms functor,
  respectively.  The existence of projective and injective resolutions ensure that derived
  functors exist, and that they can be calculated with suitable resolutions of either
  factor. As usual, for each integer $i$, we set
\[
\tor iRMN = \HH i{M\dtensor RN}\qnd \ext iRMN = \HH {-i}{\rhom RMN}
\]

The existence of resolutions allows one to attach numerical invariants to $M$.  One such
is the \emph{projective dimension} of $M$, defined to be the integer
\[
\pdim_RM = \inf\left\{n\in \BZ\left |
\begin{gathered}
  \text{there is projective complex $P$ with} \\
  \text{$P\simeq M$ and $P_i=0$ for $i\ge n+1$}
\end{gathered}\right\}\right.
\]
Thus, if $\hh M=0$, then $\pdim_RM=-\infty$, and if there is no integer $n$ for which an
appropriate $P$ exists, then $\pdim_RM=\infty$.  Similarly, $M$ has a \emph{flat
  dimension} and an \emph{injective dimension}, denoted $\fdim_RM$ and $\idim_RM$,
respectively.
\end{bfchunk}

\begin{bfchunk}{Minimal resolutions.}
\label{minimal resolutions}
Let $(R,\fm,k)$ be a local ring, with maximal ideal $\fm$ and residue field $k$. As per
the commutative algebraist's convention, `local' means also noetherian.  Each complex $M$
of $R$-modules with $\hh M$ degreewise finite and bounded below has a \emph{minimal free
  resolution}: a {\quism} $F\to M$ with
\[
F:= \cdots\to R^{b_{i+1}} \to R^{b_i}\to R^{b_{i-1}}\to\cdots
\]
where the $b_i$ are integers, and $\dd(F)\subseteq\fm F$; these properties ensure that
$b_i=0$ for $i < \inf M$. Such a resolution is unique up to isomorphism of complexes.
Furthermore, $F$ is a direct summand of any projective resolution of $M$.  In particular,
$\pdim_RM\leq n$ if and only if in the minimal free resolution $F$ of $M$ one has $F_i=0$
for $i\ge n+1$. Roberts' monograph \cite{Pr:book} is an excellent reference for the
preceding discussion.
\end{bfchunk}

In the sequel, homologically finite complexes of finite projective dimension play an
important role, and we use various tests for detecting them. Most of these are available
in the literature in a form suitable for our needs; some are not, and so are recorded
here.  The one below is immediate when we compute the derived functors in question using
the minimal resolution of $M$.

\begin{itchunk}
\label{small:localtest}
Let $R$ be a local ring and $M$ a complex of $R$-modules which is homologically finite and
bounded below. One has equalities
\[
\pdim_RM = \sup(k\dtensor RM) = - \inf\left(\,\rhom RMk\right)
\]
When $\pdim_RM$ is finite, it is also equal to $-\inf\left(\,\rhom RMR\right)$.  \qed
\end{itchunk}

One application of this test for projective dimension is:

\begin{itchunk}
\label{small:tensorproduct}
If $M$ and $N$ are complexes of $R$-modules with $\hh M$ and $\hh N$ degreewise finite and
bounded below, then
\[
\pdim_R(M\dtensor RN) = \pdim_R M + \pdim_R N
\]
\end{itchunk}

Indeed, this follows from \eqref{small:localtest}, since in $\lcat R$ one has isomorphisms
\[
k\dtensor R{(M\dtensor RN)} \simeq (k\dtensor RM)\dtensor R N \simeq (k\dtensor RM)
\dtensor k(k\dtensor RN)
\]
where the second one is contained in the remark below, which, in turn, is evident when one
uses resolutions of $M$ to realize the derived functors.

\begin{itchunk}
\label{tensork}
Let $(R,\fm,k)$ be a local ring and $M$ a complex of $R$-modules. Then $k\dtensor RM$,
$\rhom RMk$, and $\rhom RkM$ are isomorphic, in $\lcat R$, to graded $k$-vector spaces
with zero differentials.\qed
\end{itchunk}

The next result is a version of \eqref{small:localtest} for non-local rings; it is similar
to \cite[(5.3.P)]{AF:hd}.

\begin{proposition}
\label{small:globaltest}
Let $R$ be a noetherian ring and $M$ a complex of $R$-modules with $\hh M$ degreewise
finite and bounded below.  If the following conditions hold:
\begin{enumerate}[\quad\rm(a)]
\item $\rhom RM{R/\fm}$ is homologically bounded below for each maximal ideal $\fm$ of
  $R$;
\item $\rhom RMR$ is homologically bounded below.
\end{enumerate}
then $\pdim_RM=-\inf \rhom RMR$, and in particular, $\pdim_RM<\infty$.
\end{proposition}

\begin{proof}
  Set $d= - \inf \rhom RMR$; possibly, $d=-\infty$, which happens when $\rhom RMR\simeq
  0$.  Evidently, $\pdim_RM\ge d$.  For the reverse inequality, it is enough to check that
  $\pdim_{R_\fp}M_\fp\le d$ for each prime ideal $\fp$ of $R$; see \cite[(5.3.P)]{AF:hd}.
  Since projective dimension is non-increasing upon localization, it suffices to verify
  this last inequality at each maximal ideal of $R$.
  
  Fix a maximal ideal $\fm$ of $R$ and let $k$ be the residue field $R/\fm$.  Since $R$ is
  noetherian and $\hh M$ is degreewise finite and bounded below, one obtains
\begin{alignat*}{2}
  &\ext i{R_\fm}{M_\fm}k \cong \ext iRM{R/\fm}_\fm =0& &\foral i\gg 0 \\
  &\ext i{R_\fm}{M_\fm}{R_\fm} \cong \ext iRMR_\fm= 0& &\foral i\ge d+1
\end{alignat*}
The complex $M_\fm$ of $R_\fm$-modules is homologically degreewise finite and bounded
below, so the vanishing in the first row implies that $\pdim_{R_\fm}M_\fm$ is finite, by
\eqref{small:localtest}. Given this, the vanishing in the second row yields
$\pdim_{R_\fm}M_\fm\le d$, again by \eqref{small:localtest}.  This is the desired
conclusion.
\end{proof}

\begin{bfchunk}{Koszul complexes.}
  Let $\bsx$ be a finite set of elements in $R$. For each complex $M$ of $R$-modules, the
  \emph{Koszul complex} on $\bsx$ with coefficients in $M$ is denoted $\koszul \bsx M$.
  One way to construct it is by using the recipe: $\koszul{\emptyset}M = M$, and for $\bsy
  = \bsx\sqcup \{y\}$ the complex $\koszul{\bsy}M$ is the mapping cone of the endomorphism
\[
\koszul{\bsx} M \to \koszul \bsx M \where m\mapsto y m\quad \text{for $m\in \koszul\bsx
  M$}
\]
The reader may consult Bruns and Herzog~\cite[\S(1.6)]{BH} for details.

Let $(R,\fm,k)$ be a local ring and $\bsx$ a minimal generating set for $\fm$; Nakayama's
lemma implies $\card(\bsx) = \rank_k (\fm/\fm^2)$.  One refers to $\koszul\bsx R$ as
\emph{the Koszul complex of $R$}. No reference is made to $\bsx$ because if $\bsy$ is
another minimal generating set for $\fm$, then $\koszul\bsx R\cong \koszul \bsy R$ as
complexes of $R$-modules; see the discussion after \cite[(1.6.20)]{BH}.
\end{bfchunk}

The following properties of the Koszul complex are well known.

\begin{itchunk}
\label{koszul:finiteness}
Let $R$ be a ring, $M$ a complex of $R$-modules, and $\bsx\subset R$ a finite
set.
\begin{enumerate}[{\quad\rm(a)}]
\item The ideal $(\bsx)$ annihilates $\hh{\koszul \bsx M}$.
\item If $\hh M$ is degreewise noetherian, then so is $\hh{\koszul \bsx M}$.
\item If $\hh M$ is degreewise artinian, then so is $\hh{\koszul \bsx M}$.
\item If $\amp(M)$ is finite, then so is $\amp(\koszul \bsx M)$; the converse holds when
  $\bsx$ is in the Jacobson radical of $R$ and the $R$-module $\hh M$ is degreewise
noetherian.
\end{enumerate}
\end{itchunk}

Indeed, (b) and (c) are immediate from the description of $\koszul\bsx M$ as an iterated
mapping cone, as is the first claim in (d); for its converse, one uses  Nakayama's
lemma.  For (a), one constructs  an explicit homotopy between
multiplication by $x$ and by zero on $\koszul\bsx M$; see \cite[(1.6.5)]{BH}.

\section{Thick subcategories of  derived categories}
\label{section:Thicksubcategories}

Let $R$ be a ring; recall that $\lcat R$ denotes the derived category of $R$-modules.
This work makes critical use of the \emph{triangulated} structure on $\lcat R$; once
again, \cite{Ha} and \cite{Ve} are suitable references for this topic.  The aim of this
section is to recall certain important notions based on the triangulated structure, and to
provide numerous examples to illustrate the ideas involved.

\begin{chunk}
\label{thick:definition}
A \emph{thick subcategory} $\calt$ of $\lcat R$ is a non-empty full subcategory such that:
\begin{enumerate}[\quad\rm(a)]
\item $\calt$ is closed under isomorphisms in $\lcat R$;
\item In any exact triangle, if two of the objects are in $\calt$, then so is the third;
\item Every direct summand of an object in $\mathcal T$ is also in $\mathcal T$.
\end{enumerate}
Thick subcategories are synonymous with the \emph{\'epaisse subcategories} of Verdier
\cite{Ve}, although the definition in \emph{loc.~cit.} differs from the one given here;
the formulation above is due to Rickard \cite{Ri}.  If a thick subcategory $\calt$ is
closed under homotopy colimits, then it is said to be \emph{localizing}; see \cite{Ne}.
\end{chunk}

\begin{example}  
\label{thick:classes}
The following subcategories of $\lcat R$ are thick:
\begin{enumerate}[\quad\rm(i)]
\item The subcategory of complexes of finite flat dimension.
\item The subcategory of complexes of finite projective dimension.
\item The subcategory of complexes of finite injective dimension.
\end{enumerate}
When $R$ is noetherian, the following subcategories are also thick:
\begin{enumerate}[\rm(iv)]
\item[\rm(iv)] The subcategory of complexes whose homology is (degreewise) finite;
\item[\rm(v)] The subcategory of complexes whose homology is (degreewise) artinian.
\end{enumerate}

We verify that the subcategory defined in (i), which we denote $\calf$, satisfies the
conditions in \eqref{thick:definition}, and hence thick; the arguments for the other
subcategories listed above are similar. The simplest method of checking that $\calf$ is
thick is to take recourse to a characterization of finite flat dimension in terms of
derived functors: a complex $F$ of $R$-modules has flat dimension $\leq n$ if and only if
$\sup(F\dtensor RX)\leq n$ for each $R$-module $X$; see \cite[(2.4.F)]{AF:hd}.

Clearly, $\calf$ is closed under isomorphisms; this takes care of (a).  As to (b), let
\[
\exact LMN
\]
be an exact triangle in $\lcat R$. For each $R$-module $X$, the homology long exact
sequence arising from the exact triangle
\[
\exact {(L\dtensor RX)}{(M\dtensor RX)}{(N\dtensor RX)}
\]
implies that if any two of the numbers $\sup(L\dtensor RX)$, $\sup(M\dtensor RX)$, and
$\sup(N\dtensor RX)$, is $\leq n$, then the third is $\leq n+1$. Thus, if any two of $L$,
$M$, and $N$ are in $\calf$, then so is the third; this settles (b).  Finally, if $N$ is a
direct summand of $M$, then $N\dtensor RX$ is a direct summand of $M\dtensor RX$, so
$\sup(N\dtensor RX)\leq \sup(M\dtensor RX)$. Thus, $M$ in $\calf$ implies $N$ in $\calf$,
as desired.

Note that none of the subcategories (i)--(v) is localizing, for they are not closed under
homotopy colimits. Here is an example of a localizing subcategory; see, for example,
\cite[(5.3)]{DG} for a proof.
\end{example}

\begin{example}
  Let $I$ be an ideal in noetherian ring $R$. The subcategory of $\lcat R$ consisting of
  complexes $M$ where $\HH iM$ is $I$-torsion for each integer $i$, is localizing.
\end{example}

\begin{chunk}
  Let $M$ be a complex of $R$-modules.  The intersection of thick subcategories of $\lcat
  R$ containing $M$ is itself a thick subcategory. We refer to it as the thick subcategory
  \emph{generated} by $M$, and denote it $\rthick RM$.  It is suggestive to think of
  complexes in $\rthick RM$ as being \emph{finitely built} from $M$.
  
  In the same vein, there is a smallest (with respect to inclusion) localizing subcategory
  generated by $M$, and we speak of any object in this subcategory as being \emph{built}
  from $M$. This subcategory does not arise often in this work, so we resist the
  temptation to introduce notation for it.
\end{chunk}

\begin{example}
  Let $(R,\fm,k)$ be a local ring. A complex $M$ of $R$-modules is built from $k$ if and
  only if for each integer $i$, the $R$-module $ \HH iM$ is $\fm$-torsion; it is finitely
  built from $k$ if and only the length of the $R$-module $\hh M$ is finite.
  
  Indeed, the criterion for $M$ being built from $k$ is a special case of the preceding
  example, and once again we refer the reader to \cite{DG} for details. As to the one for
  being finitely built, it is not hard to verify that the subcategory of complexes with
  finite length homology is thick; since this subcategory contains $k$, it must contain
  also $\rthick Rk$. This means that anything finitely built from $k$ has homology of
  finite length.  Conversely, let $M$ be a complex such that the length of $\hh M$ is
  finite; thus, $\hh M$ is bounded and the length of each $\HH iM$ is finite.  Since $R$
  is local, the only simple module is $k$, so $\HH iM$ admits a finite filtration with
  subquotients isomorphic to $k$; this implies that, in $\lcat R$, each $\HH iM$ is
  finitely built from $k$. Now an straightforward induction on the amplitude of $M$ yields
  that $M$ is finitely built from $k$; see \eqref{yoga:homology}.
\end{example}

The localizing subcategory generated by $R$ is all of $\lcat R$; this is another way of
saying that each complex is isomorphic to one consisting of projective modules.  The thick
subcategory generated by $R$ plays a prominent role in our investigations, so it is
convenient to introduce some language to address complexes in it.

\begin{bfchunk}{Small complexes.}
  We say that a complex of $R$-modules is \emph{small} if it is finitely built from $R$.
  The following result reconciles the definition adopted here with those discussed in the
  introduction. Our choice of definition was dictated by closeness to the spirit of this
  paper.

\begin{itchunk}
\label{small:definitions}
Let $M$ be a complex of $R$-modules. The following conditions are equivalent:
\begin{enumerate}[\quad\rm(a)]
\item $M$ is \perf;
\item $M$ is isomorphic in $\lcat R$ to a finite complex of finite projective modules;
\item $\rhom RM-$ commutes with arbitrary direct sums in $\lcat R$.
\end{enumerate}
\end{itchunk}

This result is well known; the details are given for the sake of readability. Also, we
sketch only the proof that (a) and (b) are equivalent, for that is all we use in this
article.

(a) $\iff$ (b): It is not hard to verify that the class of complexes isomorphic in $\lcat
R$ to finite complexes of finite projectives is thick, and since it contains $R$, it
contains also $\rthick RR$, that is to say, every small complex. For the converse, it
suffices to prove each complex of the form
\[
P:= 0\to P_n\to \cdots \to P_m\to 0
\]
where each $R$-module $P_i$ is finite and projective, is in $\rthick RR$. To this end we
induce on $n-m$. First note that Property (b) in \eqref{thick:definition} implies that any
finite free module is in $\rthick RR$, and hence by Property (c) of \emph{loc.~cit.}, so
is any direct summand of such a module.  Thus, any finite projective module is in $\rthick
RR$, which settles the base case $n=m$. As to the induction step, assuming $n-m\geq 1$,
note that $P$ is the mapping cone of the morphism of complexes $\dd_n^P \col \shift{n-1}
{P_n} \to P_{\les n-1}$, so in $\lcat R$ one has a distinguished triangle
\[
\exact{P_{\les n-1}}{P}{\shift{n-1}{P_n}}
\]
Now, the induction hypothesis yields that $P_{\les n-1}$ and $\shift {n-1}{P_n}$ are in $\rthick
RR$, and hence so is $P$, by Property (b) of \eqref{thick:definition}.
\end{bfchunk}

Over a noetherian ring, a homologically finite complex is isomorphic in the derived
category to a finite complex of finite projectives if and only if it has finite projective
dimension, if and only if it has finite flat dimension.  Thus, over such rings, the
equivalence (a) $\iff$ (b) of the preceding result yields the following homological
characterization of \perf ness.

\begin{itchunk}
\label{perf=small}
Let $R$ be a noetherian ring and $M$ a homologically finite complex of $R$-modules.  The
following conditions are equivalent:
\begin{enumerate}[\quad\rm(a)]
\item $M$ is small;
\item $\pdim_RM$ is finite;
\item $\fdim_RM$ is finite.\qed
\end{enumerate}
\end{itchunk}

We record some elementary remarks concerning the process of building one thing out of
another. In presenting these, it is expedient to introduce the following notation: given
complexes $M$ and $N$ of $R$-modules, we write
\[
M\rbuilds RN
\]
to indicate that $N$ is in $\rthick RM$, that is to say, that $N$ is finitely built from
$M$.  Our reason for adopting this symbolism is that, often, properties of a complex $M$
are inherited by complexes finitely built from it. We hope that the discussion to follow
will convince the reader that the language, of one complex being built from another, and
the notation, $M\builds N$, aid intuition about these notions.

\medskip

In the following remarks, $L$, $M$, and $N$ are complexes of $R$-modules.

\begin{itchunk}
\label{yoga:transitivity}
If $M$ is finitely built from $L$, and $N$ is finitely built from $M$, then $N$ is
finitely built from $L$. In pictures:
\[
L\rbuilds R M \quad\text{and}\quad M\rbuilds RN \quad\text{entails}\quad L\rbuilds RN
\]
\end{itchunk}

The assertion above is verified by a direct appeal to the definition of $\rthick R-$.

\begin{itchunk}
\label{yoga:homology}
If\, $\hh N$ is bounded and for each integer $i$ the $R$-module $\HH iN$ is finitely built
from $M$, then $N$ is finitely built from $M$.
\end{itchunk}

Indeed, this claim is evidently true when $\amp N=0$, for then $N\simeq \HH iN$ for some
$i\in\BZ$. This is the basis of an induction on $\amp N$. As to the induction step: for
$s=\sup N$, by suitably truncating $N$ one can construct an exact triangle
\[
\shift s{\HH sN} \to N \to \wt N \to \shift {s+1}{\HH sN}
\]
where the induced map $\HH iN\to \HH i{\wt N}$ is bijective for $i\ne s$, and $\HH s{\wt
  N}=0$. Therefore, the induction hypothesis yields that $\rthick RM$ contains both
$\shift s{\HH sN}$, because its amplitude is zero, and $\wt N$, because its amplitude is
strictly less than the amplitude of $N$. Since $\rthick RM$ is closed under exact
triangles, one concludes that it contains $N$ as well.
  
Caveat: the converse to \eqref{yoga:homology} need not hold, as the next example
demonstrates.

\begin{example}
  Let $(R,\fm,k)$ be a local ring and $K$ the Koszul complex of $R$.  Given
  \eqref{perf=small}, it is clear that $K$ is finitely built from $R$, while $\HH 0K$,
  that is to say, $k$, is finitely built from $R$ if and only if the ring $R$ is regular.
\end{example}

\begin{itchunk}
\label{yoga:basechange}
Let $\vf\col R\to S$ be a homomorphism of rings, and $P$ a complex of $S$-modules.
\[
M\rbuilds R N \qquad \text{entails} \qquad \left\{
\begin{gathered}
  (P\dtensor RM)\rbuilds S(P\dtensor RN) \\
  {\rhom RPM} \rbuilds  S{\rhom RPN} \\
  {\rhom RMP} \rbuilds S{\rhom RNP}
\end{gathered} \right. 
\]
\end{itchunk}

The claims above are immediate, for the functors $P\dtensor R-$, $\rhom RP-$, and $\rhom
R-P$ from $\lcat R$ to $\lcat S$ respect the triangulated structures.

For the next remark, which is also easy to verify, recall that given a homomorphism
$\psi\col Q\to R$ of rings, one may view $M$ and $N$ as complexes of $Q$-modules, by
restriction along $\psi$.

\begin{itchunk}
\label{yoga:restriction}
Let $\psi\col Q\to R$ be a homomorphism of rings.  If, as complexes of $R$-modules, $N$ is
finitely built from $M$, then the same is true when they are viewed as complexes of
$Q$-modules. In pictures:
\[
M\rbuilds R N \qquad\text{implies}\qquad M\rbuilds QN
\]
\end{itchunk}

Note that the converse of the preceding remark does not hold.

\begin{example}
  Let $Q$ be a ring, $\psi$ the canonical inclusion $Q\to Q[x]/(x^2)=R$, and let $N$ be
  the $R$-module $R/(x)\cong Q$.  As a $Q$-module, $R\cong N\oplus N$, so $N$ is finitely
  built from $R$.  However, it is easy to prove that as $R$-modules $N$ is not finitely
  built from $R$.
\end{example}

Next, we recall a notion of support for complexes, introduced by Foxby~\cite {Fo}.

\begin{bfchunk}{Support.}
  Let $R$ be a ring. Given a prime ideal $\fp$, we denote $\res \fp$ its residue class
  field $R_\fp/\fp R_\fp$ at $\fp$. The \emph{small support} of a complex $M$ of
  $R$-modules is the subset of $\spec R$ described by
\[
\supp_R M =\{\fp\in\spec R \mid \big(\res \fp\dtensor RM\big)\not\simeq 0\}
\]   
Usually, we drop the subscript `R' from $\supp_RM$, unless there is cause for confusion.

The small support of $M$ is a subset of the \emph{support} of $M$, which is the set
\[
\Supp_RM = \{\fp \in\spec R\mid M_\fp \not\simeq 0\}\,,
\]
whence the name `small support'.  This inclusion can be strict; for example, the small
support of $\res\fp$ is $\{\fp\}$, while its support is $\V(\fp)$. Recall that for any
ideal $I$ in $R$, the set $\{\fq\in\spec R\mid \fq\supseteq I\}$ is denoted $\V(I)$.
 
For homologically finite complexes over noetherian rings, the two notions coincide. This
is the gist of Part (d) of the remark below, which lists also other basic properties of
the small support.

\begin{itchunk}
\label{supp:properties}
Let $R$ be a ring and let $M$ and $N$ be complexes of $R$-modules.
\begin{enumerate}[\quad\rm(a)]
\item For each finite subset $\bsx$ of $R$ one has $\supp_R\koszul\bsx R=\V(\bsx)$.
\item If $N$ is built from $M$, then $\supp_R N\subseteq \supp_R M$.
\item $\supp_R (M\dtensor RN) = \supp_R M \cap \supp_R N$.
\item When $R$ is noetherian and $\hh M$ finite, one has equalities
\[
\supp_R M = \bigcup_{i\in \BZ}\Supp_R \HH iM = \Supp_R M
\]
In particular, $\supp_R M$ is a closed subset of $\spec R$.
\end{enumerate}
\end{itchunk}

Indeed, (a) is a routine verification, as is (c), as long as we keep \eqref{tensork} in
mind.  For (d), see the discussion preceding \cite[(2.6)]{Fo}, while the key to the proof
of (b) is that (derived) tensor products commute with homotopy colimits.
\end{bfchunk}

We end this section by recalling the following remarkable theorem discovered by M.~Hopkins
\cite[(11)]{Ho}; also see Neeman \cite[(1.2), (2.8)]{Ne} and Thomason \cite[(3.15)]{Th}.
It is a powerful test for detecting when one complex is (finitely) built from another.

\begin{itchunk}
\label{hopkins}
Let $R$ be a noetherian ring and let $M$ and $N$ be complexes of $R$-modules.  If $\supp
N\subseteq \supp M$, then $N$ is built from $M$; it is also finitely built from $M$ when
both $M$ and $N$ are \perf. \qed
 \end{itchunk}

\section{Virtually small complexes}
\label{section:Virtuallysmallcomplexes}
Let $R$ be a commutative ring. This section introduces two new classes of complexes, which
are the focus of study in this article, and develops some of the elementary theory
concerning them.

\begin{bfchunk}{Virtually small complexes.}
\label{vsmall:defn}
A complex $M$ of $R$-modules, with $\hh M\ne 0$, is \emph{\vsmall} if $\rthick RM$
contains a non-zero small complex; in other words, there is a complex $W$ of $R$-modules,
with $\hh W\ne 0$, such that
\begin{enumerate}[{\quad\rm(i)}]
\item $W$ is finitely built from $R$, and
\item $W$ is finitely built from $M$.
\end{enumerate}
The complex $M$ is \emph{\psmall} if in addition to (i) and (ii) one has that
\begin{enumerate}[{\rm(iii)}]
\item $M$ is built from $W$
\end{enumerate}
By decree, a homologically zero complex is \psmall, and hence also \vsmall.  The complex
$W$ is said to be \emph{a witness that $M$ is \vsmall}, or \emph{a witness that $M$ is
  \psmall}, as is appropriate. When the property in question is clear from the context, we
state simply that \emph{$W$ is a witness for $M$}.
\end{bfchunk}

Proxy small complexes were introduced in \cite{DGI}, where they played a crucial role in
the development of a framework for studying duality phenomenon in algebra and topology.
From the perspective of homological algebra, the notion of `virtual smallness' is a
natural generalization of that of `smallness', for the behaviour of \vsmall\, complexes is
similar to that of \perf\, complexes; this is borne out by the results established in this
article. First, the obvious example.

\begin{itchunk}
\label{small=psmall}
Each \perf\, complex is \psmall, and hence also \vsmall. \qed
\end{itchunk}

In particular, when the ring $R$ is regular, every homologically finite complex is
\psmall. However, there are many more; confer Corollary \eqref{psmall:fdim} and
Proposition \eqref{vsmall:regular}.  Even over an arbitrary ring there is a plethora of
\psmall\, complexes.

\begin{example}
\label{psmall:sums}
Let $X$ be complex of $R$-modules and set $M=R\oplus X$. Since $R$ is a direct summand of
$M$, it is contained in $\rthick RM$. Hence $M$ is \vsmall; it is even \psmall\, because
$M$ is built from $R$. Note that $M$ is \perf\, if and only if $X$ is \perf.
\end{example}

One of the difficulties in checking whether a complex is \vsmall\, or \psmall\, is that
there are no canonical witnesses. However, the situation is under better control over
noetherian rings; this is the content of the next two propositions.

\begin{proposition}
\label{psmall:witness}
Let $R$ be a noetherian ring and $M$ a complex of $R$-modules. The following conditions
are equivalent.
\begin{enumerate}[{\quad\rm(a)}]
\item The complex $M$ is \psmall.
\item The subset $\supp M\subseteq \spec R$ is closed, and each \perf\, complex $P$ with
  $\supp P\subseteq \supp M$ is finitely built from $M$.
\item There is a finite subset $\bsx$ of $R$ with $\V(\bsx)=\supp M$ and $\koszul\bsx R$
  finitely built from $M$.
\end{enumerate}
\end{proposition}
\begin{proof}
  (a) $\implies$ (b): Let $W$ be a witness that $M$ is \psmall. Thus, $W$ and $M$ are
  built from each other, so $\supp W = \supp M$, by (\ref{supp:properties}.b). But then,
  $W$ is homologically finite so (\ref{supp:properties}.d) yields that $\supp W$ is a
  closed subset of $\spec R$, and hence the same is true of $\supp M$.  Since the
  complexes $W$ and $P$ are \perf\, and $\supp P\subseteq \supp M$, one obtains from
  \eqref{hopkins} that $P$ is finitely built from $W$. Therefore,
  \eqref{yoga:transitivity} implies that $P$ is finitely built also from $M$, as claimed.
  
  (b) $\implies$ (c): Take $\bsx$ to be a finite set with $\V(\bsx)=\supp M$.
  
  (c) $\implies$ (a): Given that $\koszul\bsx R$ is \perf\, it suffices to verify that $M$
  is built from it. By (\ref{supp:properties}.a), the support of $\koszul\bsx R$ is
  $\V(\bsx)$, which equals $\supp M$, so (\ref{hopkins}) yields the desired conclusion.
\end{proof}

There is a recognition principle also for virtual smallness; it is especially useful over
local rings.

\begin{proposition}
\label{vsmall:witness}
Let $R$ be a noetherian ring and $M$ a complex of $R$-modules with $\hh M\ne 0$.  The
following conditions are equivalent.
\begin{enumerate}[{\quad\rm(a)}]
\item The complex $M$ is \vsmall.
\item The subset $\supp M\subseteq \spec R$ contains a closed point $\{\fm\}$ with the
  property that each \perf\, complex supported on $\{\fm\}$ is finitely built from $M$.
\item There is a finite set $\bsx\subset R$ such that $(\bsx)\ne R$ and $\koszul\bsx R$ is
  finitely built from $M$.
\end{enumerate}
In particular, when $R$ is local, $M$ is \vsmall\, if and only if the Koszul complex of
$R$ is finitely built from $M$.
\end{proposition}
\begin{proof}
  (a) $\implies$ (b): For each witness $W$ that $M$ is \vsmall, the subset $\supp W\subset
  \spec R$ is non-empty and closed; this follows from (\ref{supp:properties}.d). Let
  $\{\fm\}$ be a closed point in $\supp W$. If $P$ is a \perf\, complex with $\supp P
  =\{\fm\}$, it follows from \eqref{hopkins} that $P$ is finitely built from $W$, and
  hence that it is finitely built from $M$, by \eqref{yoga:transitivity}.
  
  For (b) $\implies$ (c) set $\bsx$ to be a finite generating set for the ideal $\fm$,
  while (c) $\implies$ (a) is almost by definition; all one needs to note is that $\supp
  \koszul \bsx M$ is non-empty, by (\ref{supp:properties}.a).
\end{proof}

\begin{remark}
  The preceding propositions are based on Hopkins' theorem \eqref{hopkins}, and it is for
  this reason that the bulk of the results presented in this article are restricted to
  commutative noetherian rings. There are extensions of the Hopkins' theorem to other
  categories: modules over commutative, but not necessarily noetherian, rings, sheaves on
  a scheme, and modules over group rings.  In all these the role of thick subcategories is
  played by the less restrictive class of tensor closed thick subcategories, which
  suggests that the `correct' definition of virtual smallness and \psmall ness should
  involve such subcategories.
\end{remark}

The examples given so far of \psmall, or \vsmall , complexes are based on \perf\,
complexes. Next we identify diverse classes of \vsmall\, complexes; these demonstrate that
there is a lot more to these notions than \perf ness.

The result below will be subsumed in Proposition \eqref{psmall:regtest}.  However, it
seemed worthwhile to give `ad hoc' arguments for this one, for the strategy employed is
useful also in other cases.

\begin{proposition}
\label{psmall:k}
For each local ring $(R,\fm,k)$, the $R$-module $k$ is \psmall.
\end{proposition}

\begin{proof}
  We verify that the Koszul complex $K$ of $R$ is the witness we seek, as predicted by
  Proposition \eqref{psmall:witness}. In \eqref{vsmall:defn}, condition (i) asks that $K$
  be \perf, and that we know, while (ii) is a consequence of \eqref{yoga:homology}: since
  $k$ is a field and $\hh K$ is a finite dimensional $k$-vector space, $\hh K$ is finitely
  built from $k$ in $\lcat k$, and hence also in $\lcat R$, by \eqref{yoga:restriction}.
  One way to check (iii) is to invoke (\ref{hopkins}), as $\supp k = \supp K$. Another is
  to note that $K$ is a Differential Graded $R$-algebra and $k$ a DG $K$-module, via the
  augmentation $K\to k$.  Thus, $k$ is built from $K$ as a DG $K$-module, and hence also
  as an $R$-module, by a DG analogue of \eqref{yoga:restriction}.
\end{proof}

The proposition above  takes on greater significance in the light of the following
example that shows that not everything is \vsmall.  We defer our justification to
\eqref{proof:nonexample}, for it is based on certain general results on \vsmall\,
complexes; we do not know how to verify this directly.

\begin{example}
\label{nonexample:retracts}
Set $R=k[x,y]/(x^2,xy,y^2)$. The $R$-module $R/(y)$ is \emph{not} \vsmall.
\end{example}

One fall-out of this example and the proposition above is that $k$ cannot be used as a
`test' object for detecting rings over which every homologically finite complex is
\vsmall.

Over a local ring, if a complex of finite flat dimension is homologically finite, then it
is small--see \eqref{perf=small}--and hence \vsmall. The following result extends this
remark, and illustrates the flexibility afforded by the notion of virtual smallness.  Its
proof is rather long and involved; as is explained in it, the difficulty is in dealing
with the case when the field extension $k\to l$ is infinite.  One comment: not every flat
module is \vsmall; see Example \eqref{nonexample:ffd}.

\begin{theorem}
\label{vsmall:almostfinite}
Let $\vf \col (R,\fm,k)\to (S,\fn,l)$ be a local homomorphism and $N$ a complex of
$S$-modules with $\hh N$ either finite or artinian.  If $\fdim_RN$ is finite, then $N$ is
\vsmall\, over $R$.
\end{theorem}
 
\begin{proof}  
  We may assume that $\hh N$ is an $l$-vector space of finite rank. Indeed, let $K$ be the
  Koszul complex on $S$.  By \eqref{yoga:basechange}, since $K$ is finitely built from
  $S$, the complex $K\otimes_SN$ is finitely built from $N$. This holds in $\lcat S$, and
  hence also in $\lcat R$, by \eqref{yoga:restriction}.  Consequently, it is enough to
  prove that $K\otimes_SN$ is \vsmall\, over $R$. Moreover, $\fdim_R(K\otimes_SN)$ is
  finite; one way to see this is use (\ref{thick:classes}.i), since $K\otimes_SN$ is
  finitely built from $N$. The $S$-module $\hh{K\otimes_SN}$ is either finite or artinian,
  because $\hh N$ has this property, and annihilated by $\fn$; both these deductions use
  \eqref{koszul:finiteness}. Thus, replacing $N$ with $K\otimes_SN$ accomplishes the
  desired reduction.
    
  At this point, in the special case where $l$ has finite rank over $k$, one obtains that
  $\hh N$ is finite as a $k$-vector space, and hence also as an $R$-module.  Therefore,
  the finiteness of $\fdim_RN$ yields that the complex of $R$-modules $N$ is \perf, by
  \eqref{perf=small}, and hence also \vsmall.  
  
  When the extension of residue fields is not finite the argument is a lot more involved,
  and is a classic devissage argument.  First, we reduce to the case where both $R$ and
  $S$ are complete as follows: let $\wh R$ denote the $\fm$-adic completion of $R$, and
  let $\wh S$ denote the $\fn$-adic completion of $S$. We claim that when $\wh
  S\otimes_SN$ is \vsmall\, over $\wh R$, the complex $N$ is \vsmall\, over $R$.  To see
  this, let $K$ be the Koszul complex of $R$.  Thus $\wh R\otimes_R K$ is the Koszul
  complex of $\wh R$, so when $\wh S\otimes_RN$ is \vsmall\, over $\wh R$, one has that
\[
\xymatrixrowsep{2pc} \xymatrixcolsep{2pc} \xymatrix{ N \ar@{->}[r]^-{\simeq}_-{S} & \wh
  S\otimes_S N \ar@{=>}[r]_ -{\wh R} &\wh R\otimes_R K \ar@{<-}[r]^-{\simeq}_-R & K}
\]
In this diagram, the arrow pointing right is a \quism\, because $S\to \wh S$ is flat and
$\fn\hh N=0$, while the arrow pointing left is a \quism\, because $R\to \wh R$ is flat
and $\fm \hh K=0$. The implication is given by Proposition \eqref{vsmall:witness}. In
particular, using \eqref{yoga:restriction}, one obtains that $K$ is finitely built from
$N$, that is to say, $N$ is \vsmall\, over $R$.  To complete the reduction to the complete
case, one has to verify that the flat dimension of $\wh S\otimes_SN$ over $\wh R$ is
finite. Now, since $\fdim_RN$ is finite, base change yields that $\fdim_{\wh R}(\wh
R\otimes_RN)$ is finite. On the other hand, the \quism\, $N\to \wh S\otimes_SN$ factors as
\[
N\to {\wh R \otimes_R N} \to \wh S\otimes_S N
\]
Since $\fm\hh N=0$ the first map in the diagram above is a \quism\, and hence so
is the second.  Thus, $\fdim_{\wh R}(\wh S\otimes_SN)$ is finite, as desired.

Summing up, it remains to prove: given a local homomorphism $\vf \col R\to (S,\fn,l)$ of
complete local rings and a homologically finite complex $N$ of $S$-modules, if $\fdim_RN$ is
finite, then $N$ is \vsmall\, over $R$.  We no longer assume that $\fn\cdot \hh N=0$, for
that plays no role in the sequel.

The next step is to construct a convenient factorization of $\vf$, and reduce the problem
to one about complexes over regular local rings. The local ring $R$ being
complete, Cohen's structure theorem gives a surjective map $\pi\col (P,\fp,k)\to R$ with
$P$ a regular local ring. We embed this into a commutative diagram of local homomorphisms

\[
\xymatrixrowsep{2pc}
\xymatrixcolsep{3pc} 
\xymatrix { 
    &Q\ar@{->>}[dr]^-{\pi'}\ar@/^1.5pc/[drr]^{\rho}  \\
P\ar@{->}[ur]^\psi \ar@{->}[dr]_\pi &  
        &(R\otimes_PQ)\ar@{->}[r]^\tau &S \\
    &R\ar@{->>}[ur]_-{\psi'}\ar@/_1.5pc/[urr]_{\vf}}
\]
with the following properties:
\begin{enumerate}[\quad\rm(a)]
\item $Q$ is a regular local ring and $\rho$ is a surjective homomorphism;
\item $\psi$ is a flat local homomorphism whose fibre $Q/\fp Q$ is regular;
\item $\psi'=(R\otimes_P\psi)$, $\pi'=(\pi\otimes_P Q)$, and $\tau(r\otimes q)=\vf(r)\rho(q)$.
\end{enumerate}
This diagram is constructed as follows: since $S$ is complete applying the Cohen
factorization theorem of Avramov, Foxby, and B.~Herzog~\cite[(1.1)]{AFH} to $\vf\circ\pi$
provides the surjective homomorphism $\rho$, and the flat local homomorphism $\psi$ with
regular fibre. Since $P$ is regular, the last two properties entail $Q$ regular; see
\cite[(2.2.12)]{BH}. Property (c) prescribes the remaining homomorphisms in the diagram.

Now we return to the complex of $S$-modules $N$, which you will recall is homologically
finite with $\fdim_RN$ finite; we may also assume that $\hh N\ne 0$. The desired
conclusion is that $N$ is \vsmall\, over $R$, and for this it suffices to prove: the
Koszul complex $K$ of $Q$ is \vsmall\, over $P$.

Indeed, we claim more precisely that if $W$ is a witness that $K$ is \vsmall\, over $P$,
then the complex of $R$-modules $R\otimes_PW$ is a witness that $N$ is \vsmall\, over $R$.
By construction the homomorphism $\tau$ is surjective; in particular, $\hh N$ is finite
also over the local ring $R'=R\otimes_PQ$.  The homomorphism $\psi'$ is obtained by base
change of the flat homomorphism $\psi$ along $\pi$, so it is also flat, and its fibre is
isomorphic to that of $\psi$, and hence a regular local ring.  Thus, $\tau\circ\psi'$ is a
Cohen factorization of $\vf$, so $\fdim_RN$ finite implies $\fdim_{R'}N$ finite, by
\cite[(3.2)]{AFH}. Given this, \eqref{perf=small} yields that when viewed as a complex of
$R'$-modules $N$ is \perf. Now, the complex $R'\otimes_QK$ is \perf, by base change, and
its support is the maximal ideal of $R'$, by (\ref{supp:properties}.a).  Since
$\supp_{R'}(R'\otimes_QK)\subseteq \supp_{R'}N$, and both complexes in question are \perf,
from \eqref{hopkins} one obtains the first of the following implications
 \[
 N\rbuilds {R'} (R'\otimes_Q K) = (R\otimes_PQ)\otimes_Q K \cong (R\otimes_P K)\rbuilds R
 (R\otimes_PW)
\]
The second implication holds by \eqref{yoga:basechange} because $K$ builds $W$ in $\lcat
P$; this is where the virtual smallness of $K$ over $P$ comes in.  Thus, by restriction,
$N$ builds $R\otimes_PW$ in $\lcat R$. It remains to note that the complex of $R$-modules
$R\otimes_PW$ is \perf, by base change, and that it is non-trivial, for example, by
(\ref{supp:properties}.c). This settles our claim.

To complete the proof, it remains to check that the Koszul complex $K$ of $Q$ is
\vsmall\, over $P$.  However, $Q$ is regular, by construction, so $K\simeq l$ in $\lcat
Q$, and hence also in $\lcat P$. The action of $P$ on $l$ factors through
$k$, and $k$ is a direct summand of $l$ as a $k$-vector space, and so also as an
$P$-module. Thus, $l$ finitely builds $k$ in $\lcat P$; but then, $P$ being regular,
$k$ is \perf.
\end{proof}

Here is one corollary. Note that it uses only the special case of the preceding theorem when the residue field extension is finite, so one could also provide a more concise and direct proof.

\begin{corollary}
\label{psmall:fdim}
Let $R$ be a local ring and $M$ a complex of $R$-modules with $\hh M$ either finite or
artinian. When $\fdim_RM$ is finite, in particular, when $R$ is regular, $M$ is \psmall.
\end{corollary}

\begin{proof}
  One may assume that $\hh M\ne 0$.  When $\hh M$ is finite, the finiteness of $\fdim_RM$
  yields that $M$ is \perf, by \eqref{perf=small}, and hence also \psmall. When $\hh M$ is
  artinian and $\fdim_RM$ finite, $M$ is \vsmall, by the special case $\vf=\id^R$ of
  Theorem \eqref{vsmall:almostfinite}, and then its witness is the Koszul complex $K$ of
  $R$; see Proposition \eqref{vsmall:witness}. To deduce that $M$ is \psmall, it suffices
  to prove that $M$ is built from $K$. One way to see this is to invoke \eqref{hopkins}
  since $\supp K=\{\fm\}=\supp M$, where $\fm$ is the maximal ideal of $R$; the first
  equality is by (\ref{supp:properties}.a), and the second is obtained by a direct
  calculation, keeping in mind that $\hh M$ is artinian.
\end{proof}

Over regular local rings, and for virtual smallness, one has a more drastic result.

\begin{proposition}
\label{vsmall:regular}
Let $(R,\fm,k)$ be a regular local ring. A complex of $R$-modules $M$ is \vsmall\, if and only if 
$\fm\in \supp M$.
\end{proposition}

\begin{proof}
Indeed, when $M$ is \vsmall, Proposition \eqref{vsmall:witness} yields that $\fm$ is in
$\supp M$. The converse is contained in the diagram
\[
M\rbuilds R(M\dtensor Rk) \rbuilds R k
\]
where the first implication is by \eqref{yoga:basechange}, since $k$ is finitely built
from $R$. As to the second, \eqref{tensork} yields that $M\dtensor Rk$ is isomorphic in
$\lcat R$ to a 
graded $k$-vector space, while our hypothesis ensures that this graded $k$-vector space is
non-zero, and in particular, that $k$ is a direct summand.
\end{proof}

\begin{remark}
\label{vsmall:complete0}
Another corollary to Theorem \eqref{vsmall:almostfinite} is that any \perf\, complex of
$\wh R$-modules, where $\wh R$ is the $\fm$-adic completion of $R$, is \vsmall\, over $R$.
This is subsumed in Corollary \eqref{vsmall:complete} that asserts: any complex of $\wh
R$-modules which is \vsmall\, over $\wh R$ is \vsmall\, over $R$.
\end{remark}

Next we focus on differences between \psmall ness and virtual smallness.

\begin{proposition}
  \label{psmall:complete}
  The $R$-module $\wh R$ is \psmall\, if and only if $R$ is complete.
\end{proposition}
\begin{proof}
  The non-trivial implication is that when $\wh R$ is \psmall, $R$ is complete.  Viewed as
  an $R$-module, the support of $\wh R$ is $\spec R$; this is because it is faithfully
  flat.  Thus, Proposition \eqref{psmall:witness} implies that $R$ must be a witness to
  the \psmall ness of $\wh R$, that is to say, $R$ belongs to $\rthick R{\wh R}$. Now
  consider $\Lho\fm{-}$, the left derived functor of the $\fm$-adic completion functor;
  for details, see Greenlees-May~\cite{GM} or Lipman~\cite{Li}, and let $\eps\col\id \to
  \Lho\fm{-}$ the corresponding natural transformation.  It is not hard to verify that the
  full subcategory of $\lcat R$ with objects
\[
\comp{\fm}R = \{M\in \lcat R \mid \eps(M)\col M\to \Lho\fm M\, \text{is an isomorphism}\}
\]
is thick; this is the class of $\fm$-adically complete complexes of $R$-modules.
Therefore, $\wh R$ being complete, the thick subcategory it generates is contained in
$\comp{\fm}R$. In particular, $R\in \comp{\fm}R$, that is to say, $R$ is complete.
\end{proof}
 
The preceding result, in conjunction with the Remark \eqref{vsmall:complete0}, shows that
the class of \vsmall\, complexes is considerable larger than the class of \psmall\,
complexes.  On the flip side, the latter class is stable, as the former class is not,
under important functors on the derived category, as is evidenced by the next result.
In it, $\Lch IM$ denotes the local cohomology of $M$ supported on $I$; see \cite{Li}.

\begin{theorem}
\label{psmall:completions}
Let $I$ be an ideal in noetherian ring $R$. When a complex $M$ of $R$-modules is \psmall,
so are $\Lch IM$ and $\Lho IM$.
 \end{theorem}

\begin{proof}
  We treat the case of local cohomology; completions can be treated along the same lines.
  Let $W$ be a witness that $M$ is \psmall, and let $K$ be the Koszul complex on a finite
  set of generators for the ideal $I$. We claim that $K\dtensor RW$ is a witness that
  $\Lch IM$ is \psmall.
  
  Indeed, $K\dtensor RW$ is \perf, for example, by \eqref{yoga:basechange}, because both
  $K$ and $W$ are \perf. It is also finitely built from $\Lch IM$ as can be seen by
  following the path from the north-west to the south-east corners of the the diagram:
\[
\xymatrixrowsep{2pc} \xymatrixcolsep{2pc} \xymatrix{
  \Lch IM \ar@{->}[r]^-{\eta} \ar@{=>}[d]   & M \ar@{=>}[r] & W \\
  K \dtensor R \Lch IM \ar@{->}[r]^-{K\dtensor R\eta}_-{\simeq} & K\dtensor RM \ar@{=>}[r]
  & K\dtensor RW}
\]
where $\eta$ is the canonical morphism. The implication in the bottom row is induced by
the one in the top, while the vertical implication holds because $R\builds K$; both these
deductions use \eqref{yoga:basechange}. To see that $K\dtensor R\eta$ is an isomorphism,
note the isomorphisms
\[
K\dtensor R\Lch IM \simeq \Lch IK\dtensor RM \qnd \Lch IK\simeq K
\]
where the first can be obtained from \cite[(3.3.1)]{Li}, and the second from
\cite[(3.2.1)]{Li}.

It remains to verify that $\Lch IM$ is built from $K\dtensor RW$. Keeping in mind that $M$
and $W$ are built from each other, one has equalities
\[
\supp (K\dtensor RW) = \supp K \cap \supp W = \V(I) \cap \supp M = \supp \Lch IM
\]
where the first equality is by (\ref{supp:properties}.c), the second combines
(\ref{supp:properties}.a) and (\ref{supp:properties}.b), and the last one is by
\cite[(3.1.2)]{Li} and (\ref{supp:properties}.c). It remains to invoke \eqref{hopkins}.
\end{proof}
 
It is easy to construct examples that show that the statement on virtual smallness
corresponding to the one above does not hold.

\begin{example}
  Let $R$ be a noetherian ring and let $I$ and $J$ be proper ideals in $R$ such that $R/I$
  is not \vsmall, $R/J$ is \vsmall, and $\V(I)\cap \V(J)=\emptyset$.  The $R$-module
  $M=R/I \oplus R/J$ is \vsmall, since its direct summand $R/J$ has that property, whilst
  $\Lch IM\simeq R/I$ is not.
  
  For instance, let $k$ be a field and $R=k[x,y,z]/(x^2,xy,y^2)$. Pick a non-zero element
  $a\in k$ and consider the ideals $I=(y,z)$ and $J=(x,y,z-a)$ in $R$; evidently,
  $\V(I)\cap \V(J)=\emptyset$.  It follows from Example \eqref{nonexample:retracts} that
  $R/I$ is not \vsmall. On the other hand, $R/J\cong k$, so it is \vsmall, with witness
  the Koszul complex on $\{x,y,z-a\}$.
\end{example}

The construction above points the way to the analogue of Theorem
\eqref{psmall:completions} for virtual smallness.

\begin{theorem}
\label{vsmall:completions}
Let $I$ be an ideal in a local ring $R$.  If a complex $M$ of $R$-modules is \vsmall, then
so is $\Lch IM$; the converse holds if $\Lch IM\not\simeq 0$. The corresponding statement
for $\Lho IM$ also holds.
\end{theorem}
\begin{proof}
  Once again, we only provide details for the statement concerning local cohomology.  Let
  $K$ be the Koszul complex on a finite set of generators for $I$. The crucial point in
  the proof is that if a complex $W$ is \perf\, with $\hh W\ne 0$, then the complex
  $K\dtensor RW$ has the same properties. Indeed, the \perf ness of $K\dtensor RW$ follows
  from \eqref{yoga:basechange}. As to its non-triviality: both $\supp K$ and $\supp W$ are
  non-empty closed subsets of $\spec R$, by (\ref{supp:properties}.d), and hence contain
  the maximal ideal of $R$.  Thus, the same is true of $\supp(K\dtensor RW)$, by
  (\ref{supp:properties}.c); this implies $\hh{K\dtensor RW}\ne 0$.
  
  First we verify that if $M$ is \vsmall\, with witness $W$, then $\Lch IM$ is \vsmall\,
  with witness $K\dtensor RW$. We may assume that $M\not\simeq 0$.  Arguing as in the
  proof of the preceding result, one obtains that $K\dtensor RW$ is finitely built from
  $\Lch IM$. It remains to note that, by the discussion in the last paragraph, $K\dtensor
  RW$ is \perf\, and non-zero.
  
  Assume $\Lch IM$ is non-zero and \vsmall, with witness $W$. Then one has a diagram
\[
M \builds (K\dtensor RM) \simeq (K\dtensor R\Lch IM) \builds (K\dtensor RW)
\]
where the implications are given by \eqref{yoga:basechange}, and the equivalence is
verified as in the proof of Theorem \eqref{psmall:completions}.  Since $K\dtensor RW$ is
\perf\, and non-trivial, it is a witness that $M$ is \vsmall.
\end{proof}

To end this section, we present a flat module which is not \vsmall.

\begin{example}
\label{nonexample:ffd}
Let $R$ be a noetherian ring and $\fp$ a non-maximal prime ideal.  We claim that the flat
$R$-module $R_\fp$ is not \vsmall.

Indeed, if $R_\fp$ builds a complex $W$, then $\supp W \subseteq \supp R_\fp$.  If $W$ is
\perf, then $\supp W$ is a closed subset of $\spec R$, and hence $\supp W\not \subseteq
\supp R_\fp$. Thus, $\rthick R{R_\fp}$ contains no \perf\, complexes, except those that
are homologically zero.
 \end{example}

\section{Virtually small complexes as test objects}
\label{section:testobjects}

\setcounter{subsection}{1} In this section $(R,\fm,k)$ is a local ring.  Many ring
theoretic properties of $R$ can be detected by the finiteness of suitable homological
dimensions of some `test object'. Typically, either $R$ or $k$ plays this role. For
instance, the ring $R$ is Gorenstein precisely when it has finite injective dimension when
viewed as a module over itself, while $R$ is regular precisely when $k$ has finite flat
dimension.

The main goal of this section, and the next one, is to demonstrate that, often, any
\vsmall\, complex can be used as a test object. Specialized to the situation when the
complex is \perf\, or the residue field, some of these recover well known results in the
literature; however, even then the proofs we present are new.  In each case, the crux of
the argument is similar, so we begin describe it in a sufficiently general scenario.

\begin{chunk}
\label{lambda}
Let $\lambda$ be a numerical invariant defined on $\lcat R$.  The guiding examples are as
follows.
\begin{itemize}
\item[\quad\rm(a)] Let $X$ be a complex of $R$-modules. For each complex $M$ in $\lcat R$,
  set
\[
\lambda(M) = \amp(X\dtensor RM)
\]
\end{itemize}
In the next two examples, $\psi\col Q\to R$ is a local homomorphism; keep in mind the
special case where $\psi$ is the identity on $R$.
\begin{itemize}
\item[\quad\rm(b)] For each complex $M$ of $R$-modules, set \( \lambda(M) = \fdim_QM \)
\end{itemize}

One has also the version involving injective dimensions:

\begin{itemize}
\item[\quad\rm(c)] For each complex $M$ of $R$-modules, set \( \lambda(X) = \idim_QX \)
\end{itemize}
To a given invariant $\lambda$, we associate the full subcategory of $\lcat R$ defined by
\[
\subcat \lambda R = \{X\in \lcat R\mid \lambda(X)<\infty\}
\]
Here is a noteworthy feature of examples (a)--(c) above: \emph{$\subcat\lambda R$ is a
  thick subcategory of $\lcat R$.}  The required arguments are a simple extension of those
required for verifying the assertions in \eqref{thick:classes}, which correspond to the
special case $\psi=\id^R$.
\end{chunk}

We record the following principle that is used repeatedly in this section.  It is evident
that there is a more general version, dealing with properties of complexes other than the
finiteness of some numerical invariant. However, we have opted to state it in its simplest
form for it is adequate for our present purposes.

\begin{bfchunk}{Principle.}
\label{proofprinciple}
\emph{Let $K$ be the Koszul complex on a finite set of generators of $\fm$.  Suppose
  $\lambda$ has the following properties:
\begin{enumerate}[{\quad\rm(a)}]
\item $\subcat \lambda R$ is a thick subcategory of $\lcat R$;
\item if $\lambda(K)$ is finite, then so is $\lambda(R)$.
\end{enumerate}
If $M$ is \vsmall, with $\hh M\ne 0$, and $\lambda(M)$ is finite, then $\lambda (R)$ is
finite.}
\end{bfchunk}

Indeed, when $\lambda(M)$ is finite, $\subcat \lambda R$ contains $M$, and hence also
$\rthick RM$, by hypothesis (a). If $M$ is \vsmall, then $K$ lies in $\rthick RM$, by
\eqref{vsmall:witness}, so we conclude that $K$ lies in $\subcat\lambda R$, that is to
say, $\lambda(K)$ is finite. Now hypothesis (b) yields the desired conclusion.

The rest of this section provides diverse applications of this proof-principle.

\subsection*{Homological dimensions}
The theorem below unifies and extends numerous characterizations of Gorenstein rings. The
equivalence of (a) and (b) is \cite[(13.2)]{AIM}, while (c) $\implies$ (a), applied to
$\psi=\id^R$, contains the well known result that when $R$ has a homologically finite
complex of finite injective dimension and finite flat dimension, then it is Gorenstein;
see also Corollary \eqref{gor:echar} below.

\begin{theorem}
\label{gor:char}
Let $\psi\col Q\to R$ be a local homomorphism. The following conditions are equivalent.
\begin{enumerate}[\quad\rm(a)]
\item $Q$ is Gorenstein and $\fdim_QR$ is finite;
\item $\idim_QR$ is finite;
\item $\idim_QM$ is finite for a \vsmall\, complex $M$ over $R$, with $\hh M\neq0$.
\end{enumerate}
 \end{theorem}

\begin{proof}
  As noted before, (a) $\iff$ (b) is \cite[(13.2)]{AIM}. Taking $M=R$ settles (b)
  $\implies$ (c), so it remains to verify that (c) implies (b).
  
  For each complex $X$ of $R$-modules, set $\lambda(X)=\idim_QX$.  As noted in
  \eqref{lambda}, the subcategory $\subcat\lambda R$ of $\lcat R$ is thick. In order to
  apply Principle \eqref{proofprinciple} to the desired end, it remains to verify that
  when the Koszul complex $K$ of $R$ has finite injective dimension over $Q$, so does $R$.
  With $h$ the residue field of $Q$, in $\lcat R$ one has the associativity isomorphism
\[
\rhom QhK \simeq \rhom QhR \dtensor RK
\]
The complex $\rhom QhR$ of $R$-modules is homologically degreewise finite, so the
isomorphism above and (\ref{koszul:finiteness}.d) imply: $\sup\rhom QhK$ and $\sup\rhom
QhR$ are simultaneously finite. It remains to note that for any homologically finite
complex $X$ of $R$-modules, one has $\idim_QX=\sup\rhom QhX$, by \cite[(5.5.I)]{AF:hd}.
\end{proof}

One corollary is the following characterization Gorenstein rings. We do not claim it is
new; only that it illustrates the scope of the preceding result.

\begin{corollary}
\label{gor:echar}
A local ring $(R,\fm,k)$ is Gorenstein if and only if the flat dimension of the injective
hull of $k$ is finite.
\end{corollary}

\begin{proof}
  Let $E$ be the injective hull of $k$. When $R$ is Gorenstein, $\fdim_RE$ is finite,
  since $\idim_RE$ is finite; see \cite[(3.3.4)]{La:gd}.  Conversely, since $E$ is
  artinian, if $\fdim_RE$ is finite, then it is \vsmall, by Corollary \eqref{psmall:fdim},
  so it remains to invoke (c) $\implies$ (a) of Theorem \eqref{gor:char}.
\end{proof}

The next result is an analogue for flat dimensions of (c) $\implies$ (b) in 
\eqref{gor:char}. 

\begin{theorem}
\label{test:fdim}
Let $\psi\col Q\to R$ be a local homomorphism and $M$ a \vsmall\, complex of $R$-modules
with $\hh M\ne 0$. If $\fdim_QM$ is finite, then $\fdim_QR$ is finite.
 \end{theorem}

\begin{remark}
  Let $R\to S$ be another local homomorphism and $N$ a homologically finite complex of
  $S$-modules with $\fdim_RN$ finite. Then $N$ is \vsmall\, over $R$, by Theorem
  \eqref{vsmall:almostfinite}, so the preceding result yields: when $\fdim_QN$ is finite,
  so is $\fdim_QR$. In this way we recover \cite[(3.2)]{FI}.
\end{remark}

\begin{proof}
  The proof is similar to that of (c) $\implies$ (b) in Theorem \eqref{gor:char}; the
  details are given for the sake of completeness. For each complex $X$ of $R$-modules, set
  $\lambda(X)=\fdim_QX$; note that the subcategory $\subcat\lambda R$ of $\lcat R$,
  introduced in \eqref{lambda}, is thick. Let $K$ be the Koszul complex of $R$, and $h$
  the residue field of $Q$. One has the associativity isomorphism
\[
h\dtensor QK \simeq (h\dtensor QR) \dtensor RK
\]
The $R$-module $\hh{h\dtensor QR}$ is degreewise finite, so it follows from the
isomorphism above and (\ref{koszul:finiteness}.d) that $\sup(h\dtensor QK)$ and
$\sup(h\dtensor QR)$ are simultaneously finite. Thus, \cite[(5.5)]{AF:hd} implies that
$\fdim_QK$ and $\fdim_QR$ are simultaneously finite. Principle \eqref{proofprinciple} now
yields the desired conclusion.
\end{proof}

We record a direct corollary.

\begin{corollary}
\label{vsmall:retracts}
Let $Q\arto{\psi}R\arto{\pi}Q$ be local homomorphisms with $\pi\circ\psi=\id^Q$. If $Q$
viewed as an $R$-module via $\pi$ is \vsmall, then $\fdim_QR$ is finite. \qed
\end{corollary}

This result is useful for detecting some complexes that are \emph{not} \vsmall.

\begin{example}
\label{proof:nonexample}
Let $k$ be a field, $R=k[x,y]/(x^2,xy,y^2)$. As claimed in Example
\eqref{nonexample:retracts}, the $R$-module $R/(y)=k[x]/(x^2)$ is not \vsmall.

Here is a justification: set $Q=k[x]/(x^2)$, let $\psi\col Q\to R$ be the natural
inclusion, and $\pi\col R\to Q$ the canonical surjection; evidently, $\pi\circ\psi=\id^Q$.
It is easy to see that $\fdim_QR=\infty$; for example, since $Q$ is zero dimensional,
$\fdim_QR$ finite implies, due to the Auslander-Buchsbaum Equality, that $R$ is a free
$Q$-module, which it is not. Thus, Corollary \eqref{vsmall:retracts} yields the desired
conclusion.
\end{example}

Next we turn to asymptotic invariants over local homomorphisms.
 
\subsection*{Complexity and curvature}
Let $\psi\col (Q,\fq,h)\to (R,\fm,k)$ be a local homomorphism and $M$ a homologically
finite complex of $R$-modules. Let $\bsx\subset R$ be a minimal set generating the maximal
ideal of $R/\fq R$.

Following Avramov, Iyengar, and Miller~\cite{AIM} we introduce, for each integer $n$, the
\emph{$n$th Betti number}, respectively, the \emph{$n$th Bass number}, of $M$ over $\psi$
to be the number
\begin{align*}
  &\beta^\psi_n(M) = \rank_k\, \tor nQh{\koszul \bsx M}\,,\quad \text{respectively} \\
  &\mu_\psi^n(M) = \rank_k\, \ext nQh{\koszul \bsx M}
\end{align*}
Note that when $\psi=\id^R$, these coincide with the classical Betti numbers and Bass
numbers of $M$ over $R$. In \cite{AIM}, the reader will find an in depth study of
invariants over $\psi$, as well as of the ones below which reflect their asymptotic
behaviour.

The \emph{complexity} of $M$ over $\psi$ is the number
\begin{gather*}
  \cxy {\psi}M = \inf\left\{d \in\BN \left|
\begin{gathered}
  \text{there exists a number $c\in\BR$ such that}\\
  \beta^{\psi}_n(M)\leq c n^{d-1}\text{ for all $n\gg0$}
\end{gathered}
\right\}\right.
\end{gather*}
The \emph{curvature} of $M$ over $\psi$ is the number
 \[
 \curv {\psi}M = \limsup_n \sqrt[n]{\beta^\psi_n(M)}
\]
Replacing Betti numbers by Bass numbers one obtains the \emph{injective complexity}
$\injcxy \psi M$, and the \emph{injective curvature} $\injcurv\psi M$, of $M$ over $\vf$.

\begin{lemma}
\label{asymptotes:thick}
Let $M$ be a homologically finite complex of $R$-modules. For each complex $X$ in $\rthick
RM$, the following inequalities hold:
\[
\cxy \psi X \leq \cxy \psi M \qnd \curv \psi X \leq \curv \psi M
\]
\end{lemma}

\begin{proof}
  We give the argument for the claim about complexity; the one for curvature is similar.
  The key point is to verify that
\[
\calt = \{X \in \lcat R\mid \text{$X$ is homologically finite and $\cxy \psi X\leq
  \cxy\psi M$}\}
\]
is a thick subcategory of $\lcat R$. For then, since $M$ is in $\calt$, so would anything
in $\rthick RM$, which is the desired result.

Let $\bsx\subset R$ be a minimal set generating the maximal ideal of $R/\fq R$.

When $X\simeq Y$ it is immediate that $X$ and $Y$ are in $\calt$ simultaneously.  In
$\lcat R$, if $X$ is a direct summand of $Y$, then $h\dtensor Q\koszul\bsx X$ is a direct
summand of $h\dtensor Q\koszul\bsx Y$. Thus, passing to homology, for each integer $n$ one
obtains
\[
\rank_k\, \tor nQh{\koszul\bsx X}\leq \rank_k\, \tor nQh{\koszul\bsx Y}
\]
that is to say, $\beta^{\psi}_n(X)\leq \beta^{\psi}_n(Y)$.  This yields $\cxy\psi X\leq
\cxy\psi Y$, hence when $Y$ is in $\calt$, so is $X$.  Thus, it remains to establish that
in $\lcat R$, given an exact triangle
\[
\exact XYZ
\]
when any two of $X$, $Y$, or $Z$ are in $\calt$, so is the third. To see this, note that
applying $\koszul\bsx R\otimes_R -$ and then $h\dtensor Q -$ to the triangle above, and
passing to homology leads to a long exact sequence
\begin{gather*}
  \cdots\to \tor nQh{\koszul\bsx X} \to \tor nQh{\koszul\bsx Y} \to
  \tor nQh{\koszul\bsx Z}\to \\
  \to \tor {n-1}Qh{\koszul\bsx X}\to \cdots
\end{gather*}
of $k$-vector spaces. Thus, computing ranks yields inequalities
\begin{gather*}
  \beta^\psi_n(X) \leq  \beta^\psi_n(Y) + \beta^\psi_{n+1}(Z)\\
  \beta^\psi_n(Y) \leq  \beta^\psi_n(X) + \beta^\psi_n(Z) \\
  \beta^\psi_n(Z) \leq \beta^\psi_{n-1}(X) + \beta^\psi_n(Y)
\end{gather*}
From these it is clear that if any two of $\cxy \psi X$, $\cxy\psi Y$, or $\cxy\psi Z$ are
less than $\cxy\psi M$, then so is the third. Thus, the subcategory $\calt$ is thick, as
desired.
\end{proof}

\begin{theorem}
\label{test:asymptotes}
Let $\psi\col Q\to R$ be a local homomorphism and $M$ a \vsmall\, complex of $R$-modules,
with $\hh M$ finite.  Each homologically finite complex of $R$-modules $Y$ satisfies
inequalities:
\[
\cxy \psi Y\leq \cxy RY + \cxy {\psi}M \qnd \curv \psi Y\leq \max\{\curv RY,\curv
{\psi}M\}
\]
In particular, $\cxy\psi R \leq \cxy\psi M$ and $\curv\psi R\leq \curv\psi M$.
\end{theorem}

\begin{remarks}
  One has also the following injective analogues of the theorem above:
\[
\injcxy \psi Y\leq \cxy RY + \injcxy {\psi}M \qnd \injcurv \psi Y\leq \max\{\curv
RY,\injcurv {\psi}M\}
\]
The argument is similar to the one for complexity and curvature.

Another remark is that a homologically finite complex of $R$-modules has finite flat
dimension over $Q$ if and only if its complexity over $\psi$ is zero; see
\cite[(7.1.3.1)]{AIM}. Thus, the theorem above is a quantitative extension of Theorem
\eqref{test:fdim}, restricted however to homologically finite complexes.
\end{remarks}

\begin{proof}[Proof of Theorem \emph{\eqref{test:asymptotes}}]
  Applying \cite[(9.1.1.1)]{AIM}, with $\vf =\id^R$, yields
\[
\cxy \psi Y\leq \cxy RY + \cxy {\psi}R \qnd \curv \psi Y\leq \max\{\curv RY,\curv
{\psi}R\}
\]
Thus, it suffices to verify the desired inequalities in the special case $Y=R$.  Since $M$
is \vsmall, $\rthick RM$ contains $K$, the Koszul complex of $R$; this is by Proposition
\eqref{vsmall:witness}. Therefore, Lemma \eqref{asymptotes:thick} yields inequalities
\[
\cxy\psi K \leq \cxy \psi M \qnd \curv\psi K \leq \curv\psi M
\]
It remains to note that $\cxy\psi K=\cxy\psi R$ and $\curv\psi K=\curv\psi R$, by
\cite[(7.2.1)]{AIM}.
\end{proof}

The preceding theorem has implications for the descent of the complete intersection
property along local homomorphisms; see Theorem \eqref{ci:descent}.

\subsection*{Finiteness of amplitudes}
Let $R$ be a local ring and $M$ a \perf\, complex of $R$-modules.  Iversen \cite{Iv} has
proved that for each homologically finite complex $Y$ of $R$-modules, one has the
Amplitude Inequality
\[
\amp(Y) \leq \amp(M\dtensor RY)
\]
In fact, Iversen established that it is equivalent to the New Intersection Theorem, which
had then been verified only for rings containing fields. Since then, Paul
Roberts~\cite{Pr}, has proved that the latter result, and hence also the Amplitude
Inequality, is valid for all local rings.

It was proved in \cite{FI} that the Amplitude Inequality is valid even when $\hh Y$ is
only \emph{degreewise} finite; no boundedness assumptions are required. In the light of
Iversen's result, this amounted to proving that when $\amp(M\dtensor RY)$ is finite, so is
$\amp(Y)$. Numerous examples were given in \cite{FI}, see also \cite{IW}, to illustrate
the utility of this unbounded version of the Amplitude Inequality.

The theorem below, in conjunction with Theorem \eqref{vsmall:almostfinite}, subsumes
\cite[Theorem III]{FI}. Our approach differs from the one in \cite{FI}, and yields a
concise, more transparent, proof even in the special case where $M$ is a \perf\, complex.

\begin{theorem}
\label{amplitude}
Let $R$ be a local ring and let $M$ be a \vsmall\, complex of $R$ modules with $\hh M\ne
0$.  Let $Y$ be a complex of $R$-modules with $\hh Y$ degreewise finite.

If $\amp(M\dtensor RY)$ is finite, then $\amp(Y)$ is finite as well.
\end{theorem}

\begin{proof} 
  For each complex $X$ of $R$-modules, set $\lambda(X)=\amp(X\dtensor RY)$; thus, the
  desired result is that $\lambda(R)$ is finite. To obtain this, we check that the
  conditions in \eqref{proofprinciple} are satisfied: (a) is a routine verification, while
  (b) is contained in (\ref{koszul:finiteness}.d).
\end{proof}

\begin{remark}
  The conclusion of the theorem above cannot be strengthened to an inequality:
  $\amp(Y)\leq \amp(M\dtensor RY)$.  Indeed, set $R=k[x,y]/(xy)$, let $M=R/(y)$ and $Y$
  the complex
\[
0 \to R\arto{x}R\to 0
\]
The ring $R$ is a complete intersection and $\hh M$ is finite, so Theorem
\eqref{vsmall:ciring} yields $M$ is \vsmall.  However, $\amp(Y)=1$ whilst $\amp(M\dtensor
RY)=0$.
\end{remark}

\section{Gorenstein dimensions}
\label{section:Gdimensions}
Let $(R,\fm,k)$ be a local ring and $M$ a complex of $R$-modules with $\hh M$ degreewise
finite and bounded below.  We denote $\gdim_R M$ the \emph{Gorenstein dimension}, often
abbreviated to \emph{G-dimension}, of $M$ over $R$.  The reader may consult Christensen's
monograph~\cite{La:gd} for details; there, to be precise, in \cite[(2.3.8)]{La:gd}, one
finds also the following characterization of complexes of finite G-dimension that is most
apposite for us. It makes it clear that every \perf\, complex has finite $G$-dimension.

\begin{itchunk}
\label{gdim:defn}
A homologically finite complex $M$ has finite $G$-dimension if and only if
\begin{enumerate}[{\quad\rm(a)}]
\item $\rhom RMR$ is homologically bounded, and
\item the canonical biduality morphism is an isomorphism:
\[
R\arto{\simeq} \rhom R{\rhom RMR}R
 \]
\end{enumerate}
When this is the case, $\gdim_RM=-\inf\rhom RMR$.
\end{itchunk}

First, we present a new characterization of \perf\, complexes; this is in preparation for
Theorem \eqref{qgor}, and the hypotheses are dictated by it. Note that condition (ii) in
the result below is satisfied when $M$ is \vsmall, since a \perf\, complex has finite
G-dimension. One may thus view it as a test for detecting \perf\, complexes from among
\vsmall\, ones.

\begin{theorem}
\label{small:gtest}
Let $R$ be a local ring. Each homologically degreewise finite and bounded below complex
$M$ of $R$-modules satisfying the following conditions is \perf:
\begin{enumerate}[{\quad\rm(i)}]
\item $\rhom RMM$ is \perf;
\item $\rthick RM$ has a homologically non-zero complex of finite G-dimension.
\end{enumerate}
If, in addition, $\rhom RMM\simeq R$, then $M\simeq \shift d R$ for some integer $d$.
\end{theorem}

\begin{remarks}
  A homologically finite complex $M$ with $\rhom RMM\simeq R$ is said to be
  \emph{semi-dualizing}. For example, $R$ is semi-dualizing, as is a \emph{dualizing}
  complex for $R$: a semi-dualizing complex of finite injective dimension. The result
  above implies that a semi-dualizing complex of finite $G$-dimension is isomorphic to
  $\shift dR$, for some $d\in\BZ$; this is contained in \cite[(8.3)]{La:sdc}.
  
  Theorem \eqref{small:gtest} suggests that it would be fruitful to study those complexes
  $M$ for which $\rthick RM$ contains a complex of finite $G$-dimension.
\end{remarks}

The proof is given in \eqref{proof:gtest} below. It uses an elementary remark concerning
minimal resolutions:

\begin{chunk}
\label{series:pb}
The \emph{Poincar\'e series} and the \emph{Bass series} of $M$ are the formal Laurent
series
\begin{gather*}
  \rpoin RM = \sum_{n\in \ZZ}\rank_k \tor nRkM\cdot t^n \\
  \rbass RM = \sum_{n\in \ZZ}\rank_k \ext nRkM\cdot t^n
\end{gather*}
When there is no confusion about the ambient ring, we drop it from the subscript and the
superscript.  As noted in \eqref{minimal resolutions}, the complex $M$ is \perf\, if and
only if its minimal resolution is finite, which translates to the condition that $\poin M$
is a Laurent polynomial. In the sequel, it is useful to keep in mind the following
stronger, but no more difficult to establish, version of this remark: $\poin M =
\sum_{n\in \ZZ}b_nt^n$ if and only if $M$ is isomorphic to a complex of the form
\[
F:= \cdots\to R^{b_{n+1}} \to R^{b_n} \to R^{b_{n-1}}\to \cdots
\]
where $R^{b_n}$ sits in homological degree $n$, and $\dd(F)\subseteq \fm F$.
\end{chunk}

\begin{chunk}\textit{Proof of Theorem \emph{\eqref{small:gtest}}.}
\label{proof:gtest}
By (ii), in $\rthick RM$ there is a complex $G$ of $R$-modules with $\hh G\ne 0$ of finite
G-dimension.
  
First we prove that $G$ is \perf. To this end, note that \eqref{yoga:basechange} yields
\[
\rhom RMM \rbuilds R \rhom RMG \rbuilds R \rhom RGG
\]
Hence $\rhom RGG$ is finitely built from $\rhom RMM$, so $\rhom RGG$ is also \perf; here
one is using \eqref{yoga:transitivity}.  By definition \eqref{gdim:defn}, the biduality
morphism associated to $G$ is an isomorphism, and this leads to the first of the
isomorphisms:
\begin{align*}
  \rhom RkG &\simeq \rhom Rk{\rhom R{\rhom RGR}R} \\
  &\simeq \rhom R{\rhom RGR \dtensor Rk}R \\
  &\simeq \rhom k{\rhom RGR\dtensor Rk}{\rhom RkR}
\end{align*}
The second and the third isomorphisms are by adjunction; the last one uses also
\eqref{tensork}. Set $G^\ast = \rhom RGR$; this complex is homologically finite since the
G-dimension of $G$ is finite, so the isomorphisms above yield an equality of formal
Laurent series
\[
\bass G = \poin {G^*} \bass R\,.
\]
Set $\mathcal{E}=\rhom RGG$, and consider now the following isomorphisms
\begin{align*}
  \rhom Rk{\mathcal E} &\simeq \rhom RkR \dtensor R \mathcal{E} \\
  &\simeq \rhom RkR \dtensor k (k\dtensor R\mathcal{E})
\end{align*}
where the first one holds because $\mathcal E$ is small and the second follows from
\eqref{tensork} and the associativity of tensor products. On the other hand, adjunction
implies that
\begin{align*}
  \rhom Rk{\mathcal E} &\simeq \rhom R{G\dtensor Rk}G\\
  &\simeq \rhom k{G\dtensor Rk}{\rhom RkG}
\end{align*}
The two sets of isomorphisms above imply the equality of formal Laurent series:
\[
\bass R\poin{\mathcal{E}} = \poin G \bass G\,.
\]
Combining the displayed equalities above of formal Laurent series yields
\[
\bass R\poin{\mathcal{E}} = \poin G \poin {G^\ast}\bass R
\]
Since the ring of Laurent series is a domain, this equality is equivalent to
\[
\poin{\mathcal E}=\poin G \poin{G^\ast}\tag{$\ast$}
\]
The coefficients of all the Laurent series that appear in the equality are non-negative
integers and $\poin{\mathcal E}$ is a Laurent polynomial, because $\mathcal E$ is \perf.
Thus, $\poin G$ must also be a Laurent polynomial, that is to say, $G$ must be \perf, as
desired.

Next we verify that $M$ is \perf. We know now that the complex $G$ is \perf, so
\[
\rhom RGM \simeq \rhom RGR\dtensor RM
\] 
On the other hand, $\rhom RGM$ is finitely built from $\rhom RMM$, since $G$ is finitely
built from $M$, by hypothesis (ii). Thus, hypothesis (i) and the isomorphism above imply
that $\rhom RGR\dtensor RM$ is \perf. Since the ring $R$ is local and both $\rhom RGR$ and
$M$ are homologically degreewise finite and bounded below---the first because it is
\perf\, and the second by hypothesis---one deduces that $M$ is \perf; see
\eqref{small:tensorproduct}.

Finally, when $\rhom RMM\simeq R$ as well, in the analogue for $M$ of the equality
($\ast$) above, one has $\poin{\mathcal E}=1$, and that implies $\poin M = t^d$ for some
integer $d$. Therefore, $M\simeq \shift dR$, by \eqref{series:pb}.
 \end{chunk}

 In the remainder of this section, $\psi\col Q\to R$ is a local homomorphism.

\subsection*{G-dimension over $\psi$}
In \cite{IW}, Iyengar and Sather-Wagstaff introduce an invariant $\gdim_\psi M$, called
the \emph{G-dimension of $M$ over $\psi$}. It equals $\gdim_RM$ in the special case when
$\psi$ is the identity homomorphism of $R$.  In \cite{IW}, many classical results
on G-dimension over local rings are extended to this setting.  The result below, contained
in \cite[(5.1)]{IW}, is one such.

\begin{Theorem}
  If the complex $M$ of $R$-modules is \perf, with $\hh M\ne 0$, and $\gdim_\psi M$ is
  finite, then $\gdim_\psi R$ is finite as well.
\end{Theorem}

This theorem has been an important technical tool in the development in \cite{IW}; see,
for example, \cite[\S6]{IW}. This may explain the relevance of the next result, which is
analogous to Theorem \eqref{test:fdim}.

\begin{theorem}
\label{test:gdim}
Let $\psi\col Q\to R$ be a local homomorphism and $M$ a homologically finite complex of
$R$-modules with $\hh M\ne 0$. If $M$ is \vsmall\, and $\gdim_\psi M$ is finite, then
$\gdim_\psi R$ is finite.
\end{theorem}

\begin{proof}[Sketch of proof]
  One way to prove this result is to mimic \eqref{proofprinciple}: first, one observes
  that the class of homologically finite complexes of $R$-modules with finite G-dimension
  over $\psi$ is a thick subcategory of $\lcat R$. This can be deduced from the
  corresponding statement for the classical G-dimension, and it brings us to the case
  where $M$ equals $K$, the Koszul complex of $R$. Then one can either appeal to
  \cite[(5.1)]{IW}, or resort to a direct argument along the lines of the one used to
  prove that (c) $\implies$ (a) in Theorem \eqref{gor:char}.
\end{proof}

\subsection*{Quasi-Gorenstein homomorphisms}
Assume $\gdim_\psi(R)$ is finite, and let $\bass\psi$ denote the \emph{Bass series of
  $\psi$}. This is a formal Laurent series with non-negative integer coefficients that was
introduced in \cite{AF:fgd}, as was the notion of a \emph{quasi\--Gorenstein}
homomorphism: one for which $\bass \psi = t^d$ for some $d\in\BZ$; see
\cite[(7.4)]{AF:fgd}.

Suppose now that $\psi$ has a \emph{dualizing complex} $D$, in the sense of \cite{AF:fgd}.
It follows from \cite[(5.4),(7.8)]{AF:fgd} that $\psi$ is quasi-Gorenstein precisely when
$D\simeq\shift dR$ for some integer $d$.  This characterization of quasi-Gorenstein
homomorphisms is significantly generalized in the result below. Besides being of intrinsic
interest, it is plays an important role in the development in Section \ref{section:lci};
its proof is short only because most of the argument has been absorbed into Theorem
\eqref{small:gtest}.

\begin{theorem}
\label{qgor}
Let $\psi\col Q\to R$ be a local homomorphism with $\gdim_\psi R$ finite, and let $D$ be a
dualizing complex for $\psi$.  If\, $\rthick RD$ contains a homologically non-zero complex
of finite G-dimension, then $\psi$ is quasi-Gorenstein.
\end{theorem}

\begin{proof}
  The homothety $R\to \rhom RDD$ is an isomorphism; this is one of the defining properties
  of a dualizing complex for $\psi$; cf.~\cite[\S5]{AF:fgd}. Thus, when $\rthick RD$ has a
  homologically non-zero complex of finite G-dimension, Theorem \eqref{small:gtest} yields
  $D\simeq \shift dR$ for some $d\in\BZ$. As noted above, this entails $\psi$ is 
  quasi-Gorenstein.
 \end{proof}
 
\begin{remark}
  The case $\psi=\id^R$ of the theorem contains a well known characterization of
  Gorenstein rings: if $R$ has a dualizing complex of finite G-dimension, then $R$ is
  Gorenstein; see \cite[(3.4.12)]{La:gd}.
 \end{remark}

\section{Descent}
\label{section:Descent}

This section returns to the general study of \vsmall\, complexes.  We introduce some
terminology to facilitate further discussion. Given a homomorphism of rings $\vf\col R\to
S$, we say virtual smallness \emph{descends along $\vf$} if each complex of $S$-modules
that is \vsmall\, over $S$ is also \vsmall\, over $R$. In the same vein, one speaks of the
descent of \psmall ness along $\vf$. The use of the word `descent' is keeping in line with
the literature in commutative algebra, and may be less perplexing when one views $\vf$ as
an algebraic geometer does: in terms of the induced map
\[
\xymatrixrowsep{2pc} \xymatrix{
  \spec S \ar@{->}[d] \\
  \spec R }
\]
Most of the results in this section are motivated by the following:

\begin{problem}
  Identify conditions on $\vf$ which ensure that virtual smallness descends along $\vf$.
\end{problem}

We also consider the analogous problem for \psmall ness. Our first result is:

\begin{proposition}
\label{psmall:descent}
Let $\vf \col R\to S$ be a homomorphism of rings.  If $S$ is \perf\, as an $R$-module,
then both \psmall ness and virtual smallness descend along $\vf$.
\end{proposition}

\begin{proof}
  Let $N$ be \psmall\, complex of $S$-modules, with witness $W$.  Thus, from the
  properties (i)--(iii) of $W$ in \eqref{vsmall:defn} one obtains, by
  restriction~\eqref{yoga:restriction}, the corresponding properties of $W$ in $\lcat R$:
\begin{enumerate}[{\quad\rm(i)}]
\item $W$ is finitely built from $S$;
\item $W$ is finitely built from $N$, and
\item $N$ is built from $W$.
\end{enumerate}
By hypothesis, the $R$-module $S$ is finitely built from $R$, so that (i) and
\eqref{yoga:transitivity} imply: $W$ is finitely built from from $R$. This, in conjunction
with (ii) and (iii) above, yields that $W$ is a witness that $N$ is \psmall\, over $R$.

A similar argument settles the claim about virtual smallness.
\end{proof}

It may be worth drawing attention to the point that the preceding result does not require
the rings to be noetherian.  With that hypothesis, and the surjectivity of $\vf$, one can
do significantly better, and obtain a perfect `transitivity' statement.  Note that the
converses of the statements below also hold; this is a virtually trivial remark.

\begin{theorem}
\label{psmall:transitivity}
Let $\vf\col R\to S$ be a surjective homomorphism of noetherian rings.

If $S$ is \psmall\, over $R$, then \psmall ness descends along $\vf$.

If $S$ is \vsmall\, over $R$, then virtual smallness descends along $\vf$.
\end{theorem}

\begin{proof}
  We verify the claim on \psmall ness; the argument for virtual smallness is similar, and
  even a tad easier.
    Let $N$ be a \psmall\, complex of $S$-modules. By Proposition \eqref{psmall:witness},
  the subset $\supp_S N\subseteq \spec S$ is closed, that it to say, of the form $\V(I)$
  for some ideal $I$ of $S$. Pick a finite set $\bsx$ of elements in $R$ that map to a
  generating set for $I$, and let $K$ be the Koszul complex on $\bsx$. The complex
  $S\otimes_RK$ is the Koszul complex over $S$ on the set $\bsx S$, and hence it is
  \perf\, and its support is $\V(I)$; see (\ref{supp:properties}.a). Thus, Proposition
  \eqref{psmall:witness} implies that it is a witness that $N$ is \psmall\, in $\lcat S$.
  With $W$ a witness that $S$ is \psmall\, over $R$, we claim that the complex of
  $R$-modules $W\dtensor RK$ is a witness that $N$ is \psmall\,in $\lcat R$.
  
  Indeed, $W\dtensor RK$ is finitely built from $R$, since both $W$ and $K$ are finitely
  built from $R$; for example, by \eqref{yoga:basechange}. It is also finitely built from
  $N$ because
  \[
  N\rbuilds S{(S\dtensor RK)} \qnd (S\dtensor RK)\rbuilds R(W\dtensor RK)
\]
Here the first implication is valid because $S\otimes_RK$ is a witness for $N$ in $\lcat
S$, and the second follows from \eqref{yoga:basechange}, because $W$ is a witness for $S$
in $\lcat R$.  Finally, the tensor product formula for supports \eqref{supp:properties}
yields that in $\spec R$ one has equalities
\[
\supp_R(W\dtensor RK)=\supp_R W \cap \supp_R K = \supp_R S \cap \supp_R N=\supp_R N
\]
Thus, (\ref{hopkins}) implies that $N$ is built from $W\dtensor RK$, as desired.
\end{proof}

Over local homomorphisms, and for virtual smallness, one has the stronger result below; it
complements, and its proof is based on, Theorem \eqref{vsmall:almostfinite}.

\begin{theorem}
\label{vsmall:localdescent}
If $\vf \col R\to S$ is a local homomorphism with $\fdim_RS$ is finite, then virtual
smallness descends along $\vf$.
\end{theorem}

\begin{proof}
  Let $N$ be a \vsmall\, complex of $S$-modules, with witness $W$.  Since $W$ is finitely
  built from $N$, it suffices to verify that it is \vsmall\, over $R$. Then, $W$ being a
  \perf\, complex of $S$-modules, $\fdim_RS$ finite implies $\fdim_RW$ finite.  For the
  same reason, $W$ is homologically finite over $S$. Replacing $N$ by $W$ allows one to
  apply Theorem \eqref{vsmall:almostfinite}, which gives the desired result.
\end{proof}
 
Here is a direct corollary; confer Proposition \eqref{psmall:complete} which asserts that,
as an $R$-module, $\wh R$ is \emph{not} \psmall, unless $R$ is complete.

\begin{corollary}
\label{vsmall:complete}
Let $(R,\fm,k)$ be a local ring and $\wh R$ its $\fm$-adic completion.  Any \vsmall\,
complex of $\wh R$-modules is \vsmall\, over $R$. \qed
\end{corollary}

In the preceding results the hypotheses have been on the $R$-module structure on $S$; the
remaining ones focus on a property intrinsic to $S$.  The theorem below implies that when
$S$ is regular, any complex of $S$-modules whose homology is either finite or artinian, is
\psmall\, over $R$; it extends \cite[(3.2)]{DGI}.  We adopt the convention that a
(not-necessarily-local) ring is \emph{regular} if it is noetherian and has finite global
dimension.

\begin{theorem}
\label{psmall:regtest}
Let $\vf \col R\to S$ be a homomorphism of rings with $R$ noetherian and the $R$-module
$S$ finite. If $S$ is regular, then \psmall ness and virtual smallness descend along
$\vf$.
\end{theorem}

\begin{proof}
  First consider the case when $\vf$ is surjective. It suffices to verify that the
  $R$-module $S$ is \psmall; this is by Proposition \eqref{psmall:descent}. Let
  $I=\Ker(\vf)$ and $K$ the Koszul complex on a finite generating set for $I$. Thanks to
  Proposition \eqref{psmall:witness}, it is enough to check that $K$ is finitely built
  from $S$, in $\lcat R$.  For each integer $i$, the $R$-module $\HH iK$ is annihilated by
  $I$, and hence it is an $S$-module. Since $S$ is regular, it finitely builds $\HH iK$
  over $S$, and hence also over $R$.  Therefore, $K$ is finitely built from $S$ over $R$;
  see \eqref{yoga:homology}.
  
  Now we turn to the general case. Pick elements $s_1,\dots,s_e$ in $S$ that generate $S$
  as an $R$-module. Since $R$ is noetherian, for each integer $n$, the $R$-submodule of
  $S$ generated by the set $\{s_n,s_n^2,s_n^3,\dots\}$ is finitely generated; that is to
  say, there exists a monic polynomial $f_n(x_n)$ in the ring $R[x_n]$ such that
  $f_n(s_n)=0$ in $S$. Set
\[
\wt R = R[x_1,\dots,x_e]/(f_1(x_1),\dots,f_e(x_e))
\]
The ring homomorphism $R[x_1,\dots,x_e]\to S$ defined by $x_n\mapsto s_n$ is surjective
and factors through $\wt R$, to yield a homomorphism $\pi\col \wt R\to S$. Thus, one
obtains a diagram $R\arto{\eta} \wt R\arto{\pi} S$ of homomorphisms of rings with
$\pi\circ\eta=\vf$, where $\eta$ is the canonical inclusion. Because $\pi$ is surjective,
the already established part of the proposition yields that \psmall ness and virtual
smallness descend along it. By Proposition \eqref{psmall:descent}, these properties
descend also along $\eta$ because the $R$-module $\wt R$ is finite and free, as can be verified
easily. These give the desired conclusion.
 \end{proof}
 
 Once again, one has a better result for local homomorphisms.

\begin{theorem}
\label{vsmall:regulardescent}
Let $\vf \col (R,\fm,k)\to (S,\fn,l)$ be a local homomorphism.  If $S$ is regular, then
virtual smallness descends along $\vf$.
\end{theorem}

\begin{proof}
  Proposition \eqref{vsmall:witness} says that the Koszul complex $K$ of $S$ serves as a
  witness for virtual smallness in $\lcat S$. Thus, for the desired result, it suffices to
  prove that $K$ is \vsmall\, in $\lcat R$. The ring $S$ being regular, $K\simeq l$, so
  one needs to check that $l$ is \vsmall\, in $\lcat R$. However, as a $k$-module $l$
  finitely builds $k$; this is another way of saying that $k$ is a direct summand of $l$.
  Since the action of $R$ on $l$ factors through $k$, one deduces that $l$ finitely builds
  $k$ in $\lcat R$. It remains to recall that $k$ finitely builds the Koszul complex of
  $R$; see \eqref{psmall:k}.
\end{proof}

\section{Ascent}
\label{section:Ascent}
Let $\psi\col Q\to R$ be a local homomorphism.  We say that virtual smallness
\emph{ascends} along $\psi$ if each complex of $R$-modules that is \vsmall\, over $Q$ is
also \vsmall\, over $R$; the explanation for the choice of language is contained in the
discussion at the beginning of Section \ref{section:Descent}.

The main contribution of this section is Theorem \eqref{psmall:diagonal} below that gives
one criterion for ascent of the properties under consideration. Stating it requires that
we enlarge our universe from rings to (commutative) Differential Graded algebras, and
consider some notions from Section \ref{section:Thicksubcategories}, introduced for rings,
in this more general context.  The reader may refer to \cite{Av:barca} for basic results
on DG objects in the commutative setting.

\begin{chunk}
  Let $A$ be a DG algebra. Akin to the case of rings, the derived category of DG
  $A$-modules is a triangulated category, which we denote $\lcat A$. Thus, one has the
  notion of one DG $A$-module being finitely built from another, and most basic remarks in
  Section \ref{section:Thicksubcategories}, notably \eqref{yoga:transitivity},
  \eqref{yoga:basechange}, and \eqref{yoga:restriction}, apply to DG $A$-modules. In
  particular, Definition \eqref{vsmall:defn} extends to DG $A$-modules.
\end{chunk}

In what follows, the homological algebra of the `derived diagonal' homomorphism plays a
crucial role. The next paragraph takes the first steps towards explaining this statement;
see \eqref{diagonal:remark}.

\begin{chunk}
\label{diagonal:module}
Let $\phi\col Q[X]\arto{\simeq}R$ be a DG algebra model of $\psi$; the construction of DG
models is explained in \cite[\S(2.2)]{Av:barca}.  In the sequel, we write $Q[X,X]$ for DG
algebra $Q[X]\otimes_QQ[X]$.  As the graded algebra underlying $Q[X]$ is
graded-commutative, the category of DG $Q[X]$-bimodules coincides with the category of DG
(left) $Q[X,X]$-modules.  In this identification, a DG $Q[X]$-bimodule $U$ acquires a
structure of a DG $Q[X,X]$-module via the prescription
\[
(r_1\otimes r_2)\cdot u = \sign{r_2}u\, r_1ur_2
\]
for $(r_1\otimes r_2)\in Q[X,X]$ and $u\in U$. This remark applies, in particular, to
$Q[X]$ itself, with its canonical $Q[X]$-bimodule structure.
\end{chunk}

\begin{theorem}
\label{psmall:diagonal}
Let $\psi\col Q\to R$ be a homomorphism of rings and $Q[X]$ a DG model for $\psi$.  Assume
that $Q[X,X]$ is finitely built from $Q[X]$ in $\lcat{Q[X,X]}$. The following statements
hold.

\begin{enumerate}[{\quad\rm(a)}]
\item Virtual smallness and \psmall ness ascend along $\psi$.
\item If $Q$ is noetherian and $R$ is module finite over $Q$, then the $Q$-module $R$ is
  \perf, and virtual smallness and \psmall ness also descend along $\psi$.
\end{enumerate}
\end{theorem}

The proof of the theorem above is given in \eqref{diagonal:proof}. For now, we would like
to make a few remarks on its formulation.

\begin{remarks}
  Given that the category of DG Q[X,X]-modules is equivalent to that of DG
  $Q[X]$-bimodules, the ruling hypothesis in the theorem may stated as: $Q[X,X]$ is
  finitely built from $Q[X]$ as DG bimodules over $Q[X]$.  Moreover, the hypothesis does
  not depend on the choice of the model for $\psi$; this statement is clarified in Remark
  \eqref{diagonal:remark}. Therefore one may formulate the theorem above without taking
  recourse to any specific models; this is done in Remark \eqref{diagonal:remark}.
  
  Another noteworthy point is that there is nothing sacrosanct about commutative DG
  algebras: any suitable category of algebras that contains the commutative rings--for
  example, simplicial algebras--would do just as well. We have opted to work in that
  category for it suffices for the present need, and while the arguments we give do
  exploit the commutativity of the DG algebras, they can be modified to work in the
  desired generality; confer Remark \eqref{proof:categoric}.
\end{remarks}
  
One bonus of working with DG algebras rather than rings--which may be considered as
DG algebras concentrated in degree zero--is that in testing whether one DG module is built
from another one has the flexibility provided by the next result.

\begin{proposition}
\label{building:dginvariance}
Consider a morphism $\phi\col A\to A'$ of DG algebras, DG $A$-modules $U,V$, and DG
$A'$-modules $U',V'$. Let $f\col U\to U'$ and $g\col V\to V'$ be $\phi$-linear
morphisms of DG modules.

When $\hh\phi$, $\hh f$, and $\hh g$ are bijective, the DG $A$-module $V$ is finitely
built from $U$ if and only if the DG $A'$-module $V'$ is finitely built from $U'$.
\end{proposition}

\begin{proof}
  To begin with note that in $\lcat A$ the morphism $f$ factors as
\[
\xymatrixrowsep{3pc} \xymatrixcolsep{3pc} \xymatrix{ U=A\dtensor AU
  \ar@{->}[r]^-{\phi\dtensor AU} & A'\dtensor AU \ar@{->}[r]^-{A'\dtensor {\phi}f} &
  A'\dtensor{A'}U'= U'}
\]
where $A'\dtensor\phi f$ is a morphism of DG $A'$-modules.  Since $\hh{\phi}$ is
bijective, so is $\hh{\phi\dtensor AU}$, hence the bijectivity of $\hh f$ implies that of
$\hh{A'\dtensor\phi f}$. By the same token, $\hh{A'\dtensor\phi g}$ is bijective.  Now,
when the DG $A$-module $V$ is finitely built from $U$, the DG algebra version of
\eqref{yoga:basechange} yields that the DG $A'$-module $A'\dtensor AV$ is finitely built
from $A'\dtensor AU$. Thus, by the preceding discussion, we deduce that $V'$ is finitely
built from $U'$, as desired.  The converse is immediate from \eqref{yoga:restriction}.
\end{proof}
 
One useful consequence of the preceding result is the

\begin{corollary}
\label{psmall:dginvariance}
Let $\phi\col A\to A'$ be a morphism of DG algebras.  If $\hh\phi$ is bijective, then
virtual smallness and \psmall ness ascend and descend along $\phi$.
 \end{corollary}

\begin{proof}
  We deal with virtual smallness, and leave it to the reader to worry about \psmall ness.
  Let $U'$ be a DG $A'$-module.
  
  In $\lcat A$, every $A'$-linear morphism is also $A$-linear and $A\simeq A'$. Thus,
  if $U'$ is \vsmall\, over $A'$, then it is \vsmall\, over $A$.
  
  Conversely, if $U'$ is \vsmall\, over $A$, with witness $V$, then $V'=A'\dtensor AV$,
  viewed as a DG $A'$-module via the left hand factor of the (derived) tensor product, is
  a witness that $U'$ is \vsmall\, over $A'$.  Indeed, $\phi$ being an \quism, is the
  $\phi$-linear morphism
\[
\phi\dtensor A V \col V = A\dtensor AV \to A'\dtensor A V = V'
\]
Set $g=\phi\dtensor AV$. Now, Proposition \eqref{building:dginvariance} applied with
$f=\phi$ yields that in $\lcat{A'}$, the DG-module $V'$ is finitely built from $A'$, and
applied with $f=\id^{U'}$ yields that $V'$ is finitely built from $U'$.  This is the
desired conclusion.
\end{proof}

The proof of the Theorem \eqref{psmall:diagonal} uses a well known adjunction isomorphism:

\begin{chunk}
\label{diagonal:adjunction}
If $U$ and $V$ are DG $Q[X]$-modules, then in $\lcat{Q[X]}$ one has
\[
\rhom{Q[X,X]}{Q[X]}{\rhom QUV} \simeq \rhom {Q[X]}UV
\]
Indeed, $Q[X]$ being graded-commutative, $\Hom QUV$ is a DG $Q[X,X]$-module, with
\[
\big[(r_1\otimes_Qr_2)\cdot f\big](u) = \sign{r_2}f\, r_1f(r_2u)
\]
By MacLane~\cite[Chapter VI, (8.6)]{Mc}, there is an isomorphism of DG $Q[X]$-modules
\[
\Hom{Q[X,X]}{Q[X]}{\Hom QUV} \cong \Hom {Q[X]}UV\,;
\]
the derived version of this isomorphism is the one we seek.
\end{chunk}

We are now ready for the

\begin{chunk}\emph{Proof of Theorem \emph{\eqref{psmall:diagonal}}.}
\label{diagonal:proof}
(a) We focus on virtual smallness; the argument for \psmall ness is similar.  Recall that
$\phi\col Q[X]\to R$ is a DG algebra model for $\psi$. Let $\eta\col Q\to Q[X]$ be the
canonical inclusion; one has $\phi\circ\eta=\psi$.  Thanks to Corollary
\eqref{psmall:dginvariance}, since $\phi$ is a {\quism}, it suffices to prove that virtual
smallness ascends along $\eta$.

To this end, let $U$ be a DG $Q[X]$-module which is \vsmall\, over $Q$, and let $W$ be a
witness that this is so. We claim that $Q[X]\dtensor QW$, viewed as a DG $Q[X]$-module via
the left hand factor of the tensor product, is a witness that $U$ is \vsmall\, over
$Q[X]$.  Indeed, the complex $ Q[X]\dtensor QW$ is finitely built from $Q[X]$ since $W$ is
finitely built from $Q$. The argument that it is also finitely built from $U$ is contained
in the diagram:
\[
U = Q[X]\dtensor{Q[X]}U \rbuilds{Q[X]} Q[X,X]\dtensor{Q[X]}U \simeq Q[X]\dtensor QU
\rbuilds{Q[X]} Q[X]\dtensor QW
\]
where the equivalence is evident, once we recall that $Q[X,X]=Q[X]\otimes_QQ[X]$, and the
second implication is given by \eqref{yoga:basechange}, since $W$ is finitely built from
$U$.

The justification for the first implication is more delicate: as noted before, the
category of DG $Q[X,X]$-modules is the same as the category of DG $Q[X]$-bimodules.  Let
$V$ be a DG $Q[X,X]$-module. The derived tensor product $V\dtensor {Q[X]}U$, whose
construction uses the DG \emph{right} $Q[X]$-module structure of $V$, inherits from $V$
the structure of a DG \emph{left} $Q[X]$-module. Moreover, the assignment $V\mapsto
V\dtensor {Q[X]} U$ is an exact functor from $\lcat{Q[X,X]}$ to $\lcat{Q[X]}$.  Thus, if
$V$ is finitely built from $V'$ as DG $Q[X,X]$-modules, then $V\dtensor{Q[X]}U$ is
finitely built from $V'\dtensor{Q[X]}U$ as DG $Q[X]$-modules.  The special case $V=Q[X,X]$
and $V'=Q[X]$ of the last statement is the desired implication.

This completes the proof that virtual smallness ascends along $\psi$.

(b) It is enough to prove that $R$ is small over $Q$: the descent of virtual smallness and
\psmall ness follow from Proposition \eqref{psmall:descent}.  As to that, thanks to
Proposition \eqref{small:globaltest}, it suffices to verify that, homologically, $\rhom
QRN$ is bounded below for each $Q$-module $N\ne 0$. Since $Q[X]\simeq R$ in $\lcat Q$,
this is tantamount to verifying that
\[
\inf\big(\rhom Q{Q[X]}N)\, > -\infty
\]
To this end, note that in $\lcat{Q[X]}$ one has isomorphisms
\[
\rhom {Q[X,X]}{Q[X]}{\rhom Q{Q[X]}N} \simeq \rhom {Q[X]}{Q[X]}N \simeq N
\]
where the one on the left is given by \eqref{diagonal:adjunction}, while that on the right
needs no explanation.  In particular, the infimum of the DG module on the left is zero.
From this remark, our hypothesis that $Q[X,X]$ is finitely built from $Q[X]$ in
$\lcat{Q[X,X]}$, and the isomorphism
\[
\rhom Q{Q[X]}N \simeq \rhom {Q[X,X]}{Q[X,X]}{\rhom Q{Q[X]}N}
\]
we deduce: either $\rhom Q{Q[X]}N$ is zero, or its infimum is finite. \qed
\end{chunk}

The remainder of this section expatiates on the remarks following Theorem
\eqref{psmall:diagonal}.

\begin{chunk}
\label{diagonal:invariance}
Let $Q[Y]$ be another DG algebra model of $\psi$. There is then a morphism of DG
$Q$-algebras $\kappa\col Q[X]\to Q[Y]$ such that
\begin{enumerate}[\quad\rm(i)]
\item $\kappa$ is a {\quism};
\item $\kappa\otimes_Q \kappa \col Q[X,X]\to Q[Y,Y]$ is a {\quism};
\item $\kappa$ is $(\kappa\otimes_Q \kappa)$-linear.
\end{enumerate}

Indeed, the lifting property of models \cite[(2.1.9)]{Av:barca} yields a commutative
diagram of morphisms of DG $Q$-algebras
\[
\xymatrixrowsep{3pc} \xymatrixcolsep{3pc} \xymatrix{
  &Q[Y] \ar@{->}[d]^{\psi} \\
  Q[X] \ar@{->}[ur]^{\kappa} \ar@{->}[r]^-{\phi} &R}
\]
Since $\phi$ and $\psi$ are both {\quism}s, so is $\kappa$.  Therefore,
$\kappa\otimes_Q\kappa$ is also a {\quism}, because the DG $Q$-modules $Q[X]$ and $Q[Y]$
are semifree; see \cite[(1.3.3)]{Av:barca}.  Finally, it is elementary to verify that
$\kappa$ is $(\kappa\otimes_Q\kappa)$-linear.
\end{chunk}

\begin{remarks}
\label{diagonal:remark}
In what follows, when we speak about the DG $Q$-algebra $R\dtensor QR$, we mean the DG
$Q$-algebra $Q[X,X]$, and when we speak of $R$ as a DG $(R\dtensor QR)$-module, we mean
$Q[X]$ viewed as a DG $Q[X,X]$-module. Remark \eqref{diagonal:invariance} explains the
sense in which these are independent of the choice of the DG model for $\psi$.  In
addition, it is immediate from Proposition \eqref{building:dginvariance} and
\eqref{diagonal:invariance} that condition that $R\dtensor QR$ is finitely built from $R$
as $(R\dtensor QR)$-modules is also intrinsic to $\psi$. The upshot is that the hypothesis
of Theorem \eqref{psmall:diagonal} admits a rather more pleasing formulation.

\begin{Theorem}
  Let $\psi\col Q\to R$ be a homomorphism of rings. Assume that $R\dtensor QR$ is finitely
  built from $R$ in $\lcat{R\dtensor QR}$. The following statements hold.

\begin{enumerate}[{\quad\rm(a)}]
\item Virtual smallness and \psmall ness ascend along $\psi$.
\item If $Q$ is noetherian and $R$ is module finite over $Q$, then the $Q$-module $R$ is
  \perf, and virtual smallness and \psmall ness also descend along $\psi$.
\end{enumerate}
\end{Theorem}
\end{remarks}

The last item in this section is a more category theoretic proof of the following
extension of Theorem (\ref{psmall:diagonal}.b). The argument given may be more
illuminating than the one given above; besides, it lends itself to use in other contexts.

\begin{proposition}
\label{proof:categoric}
Let $\psi\col Q\to R$ be a homomorphism of rings and set $\diag QR = R\dtensor QR$. If $R$
is \psmall\, over $\diag QR$, with witness $W$, then $V=R\dtensor {\diag QR}W$, viewed as
a complex of $R$-modules via the right hand factor, has the following properties:
    \begin{enumerate}[{\quad\rm(i)}]
    \item In $\lcat Q$, the complex $V$ is finitely built from $Q$;
    \item In $\lcat R$, the complex $V$ is finitely built from $R$;
    \item In $\lcat R$, the complex $R$ is built from $V$.
    \end{enumerate}
    In particular, as a $Q$-module, $R$ is \psmall, with witness $V$.
\end{proposition}

\begin{proof}
  Property (ii) holds by base change, since $W$ is finitely built from $R$, in $\diag QR$.
  Similarly, since $R$ is built from $W$ in $\lcat{\diag QR}$, one obtains that $\diag QR$
  is built from $V$ in $\lcat R$.  However $R$ itself is built from $\diag QR$, where we
  view the latter as an $R$-module via the left hand factor. Thus, $R$ is built from $V$
  over $R$, which settles (iii).
 
  For (i), it suffices to prove--see \eqref{small:definitions}--that, given complexes of
  $Q$-modules $\{M_\lambda\}_{\lambda\in\Lambda}$, indexed by a set $\Lambda$, the
  canonical map
\[
\eta \col \colim{\lambda\in \Lambda} \rhom QV{M_\lambda} \to \rhom
QV{\colim{\lambda\in\Lambda}M_\lambda}
\]
is an isomorphism. For each DG $\diag QR$-module $X$, denote $\eta(X)$ the morphism
\[
\colim{\lambda\in \Lambda}\rhom {\diag QR}X{\rhom QR{M_{\lambda}}} \to \rhom {\diag
  QR}X{\rhom QR{\colim{\lambda\in\Lambda}M_\lambda}}
\]
induced by $\eta$.  Observe that $\eta(R)$ is an isomorphism; this is a consequence of the
isomorphism in \eqref{diagonal:adjunction}.  Thus, $\eta(X)$ is an isomorphism whenever
$X$ is finitely built from $R$; in particular $\eta(W)$ is an isomorphism.  However
\[
\rhom {\diag QR}W{\rhom QR{-}} \simeq \rhom Q{R\dtensor{\diag QR}W}{-}
\]
by adjunction, therefore $\eta(W)=\eta$, and hence $\eta$ is an isomorphism, as desired.
\end{proof}

\section{Complete intersection local rings}
\label{section:lci}
\setcounter{subsection}{1}

The results described in this section are all inspired by the question: \emph{what rings
  have the property that every homologically finite complex over them is \vsmall?} One of
our main conclusions is contained in Theorem \eqref{vsmall:ciring}; the following result
leads up to it.

\begin{theorem}
\label{vsmall:cimap}
Let $\psi\col Q\to R$ be an surjective homomorphism of noetherian rings. If $\Ker(\psi)$
can be generated by a regular sequence, then virtual smallness and \psmall ness ascend and
descend along $\psi$.
\end{theorem}
\begin{proof}
  Let $r_1,\dots,r_c$ be a regular sequence generating the ideal $\Ker(\psi)$, and set
  $Q_i= Q/(r_1,\cdots,r_i)$. Consider canonical surjections
\[
Q_0 = Q\arto{\psi_1} Q_1 \arto{\psi_2} \cdots \arto{\psi_{c-1}} Q_{c-1} \arto{\psi_c} Q_c=
R
\]
Evidently, it suffices to verify that virtual smallness and \psmall ness ascend and
descend along each of the $\psi_i$. Note that $\Ker(\psi_i)$ is generated by the regular
element $r_i$ in $Q_{i-1}$. Thus, we may assume that $R=Q/(r)$, with $r$ a regular
element.

The Koszul complex on $r$ is the DG algebra $Q[x\mid \dd(x)=r]$, with $x$ an (exterior)
variable in degree $1$. The complex of $Q$-modules underlying it is
\[
0\to Q\arto{r} Q\to 0
\]
Consider the canonical surjection $\phi\col Q[x]\to R$ of DG $Q$-algebra. Since $r$ is a
regular element, $\HH 1{Q[x]}=0$, and so $\phi$ is a {\quism}; that is to say, $Q[x]$ is a
DG algebra model for $\psi$.

As a complex of $Q$-modules, $Q[x]$, and hence $R$, is \perf, so the assertion on descent
is contained in Proposition \eqref{psmall:descent}.

To prove ascent, we verify the hypotheses of Theorem \eqref{psmall:diagonal}.  This
entails verifying that $Q[x,x]$ is finitely built from $Q[x]$ in $\lcat{Q[x,x]}$.  A
direct computation reveals that the $Q[x,x]$-module $Q[x]$ is cyclic, and its annihilator
is the DG ideal generated by the element $(x\otimes 1 - 1\otimes x)$, which we denote
$\Delta$.  Moreover, the annihilator in $Q[x,x]$ of the element $\Delta$ is the ideal
$(\Delta)$.  Consequently, one has an exact sequence of DG $Q[x,x]$-modules
\begin{gather*}
  0\to Q[x] \arto{\eta} Q[x,x]\arto{\eps} Q[x]\to 0 \\
\end{gather*}
where $\eta(1)=\Delta$ and $\eps$ is the canonical surjection.  Thus, as DG
$Q[x,x]$-modules, $Q[x,x]$ is finitely built from $Q[x]$.
\end{proof}

\begin{remark}
  Abstracting the crux of the preceding proof one obtains:
  
  \emph{Let $A$ be a graded-commutative DG algebra and $\phi\col A\to A[X]$ a semi-free
    extension, where the set $X$ is finite and concentrated in odd degrees. Then \vsmall
    ness and \psmall ness ascend and descend along $\phi$.}
  
  One benefit of this version of the result is that it contains the extension of Theorem
  \eqref{vsmall:cimap} where one assumes only that the Koszul complex on a minimal
  generating set for $\Ker(\psi)$ is acyclic.
\end{remark}

Now we apply the preceding theorem to complete intersection. First, a recap.

\begin{chunk}
\label{ci:recall}
Let $(R,\fm,k)$ be a local ring and denote $\wh R$ its $\fm$-adic completion. A
\emph{Cohen presentation} of $\wh R$ is a surjective local homomorphism $Q\tra R$ with the
ring $Q$ regular; one such always exists, by Cohen's Structure Theorem.

The ring $R$ is \emph{complete intersection} if the kernel of some, equivalently, any,
Cohen presentation $Q\tra \wh R$ is generated by a regular sequence.  In this case, if
$\pi\col Q\tra R$ is a surjective homomorphism with $Q$ a regular local ring, then
$\Ker(\pi)$ is generated by a regular sequence; see \cite[\S(2.3)]{BH} for details.
\end{chunk}

Here is the advertized result on complete intersection. As to its statement, we should
like to remark that no examples are known of complete intersections that are \emph{not}
quotients of regular rings.

\begin{theorem}
\label{vsmall:ciring}
Let $(R,\fm,k)$ be a complete intersection local ring, and $M$ a complex of $R$-modules.  

If $\fm\in\Supp M$, then $M$ is \vsmall.  If $R$ is a quotient of a regular local ring and
$\hh M$ is either finite or artinian, then $M$ is \psmall.
\end{theorem}

\begin{proof}
  When $R$ is a quotient of a regular local ring say $Q$, the kernel of the canonical
  surjection $Q\to R$ is generated by a regular sequence; this was noted in
  \eqref{ci:recall}.  The local ring $Q$ being regular, when $\hh M$ is either finite or
  artinian, Corollary \eqref{psmall:fdim} yields that $M$ is \psmall\, over $Q$; given
  this, Theorem \eqref{vsmall:cimap} implies $M$ is \psmall\, over $R$.
  
  Now we tackle the general case; in doing so, we may as well assume $\hh M\ne 0$.  
  Let  $K$ be the Koszul complex of $R$. 
 Since $K$ is finitely built from $R$, the complex
  $K\otimes_RM$ is finitely built from $M$. Moreover, using \eqref{supp:properties} one obtains
  \[
  \supp(K\dtensor RM) = \supp K \cap \supp M = \{\fm\}\cap \supp M  
  \]
  Thus, our hypothesis ensures that $\supp(K\dtensor RM)=\{\fm\}$, and in particular, that
  $\hh{K\dtensor RM}\ne 0$. Hence, it suffices to prove that $K\otimes_RM$ is \vsmall, and
  so, replacing $M$ by this complex, one may assume $\supp M =\{\fm\}$, and $\fm\hh M=0$,
  see (\ref{koszul:finiteness}.a), that is to say, $\hh M$ is a $k$-vector space.

  Let $\wh R$ be the $\fm$-adic completion of $R$. Now that $\fm\hh M=0$, the canonical
  homomorphism $M \to \wh R \otimes_RM$ of complexes of $R$-modules is a \quism.  Thus,
  $\hh{\wh R\otimes_RM}$ is again a $k$-vector space and $\fm \wh R$, the maximal ideal of
  $\wh R$, is contained in $\supp_{\wh R}(\wh R\otimes_RM)$.  Cohen's structure theorem
  gives a surjective local homomorphism $Q\to \wh R$ with $Q$ regular; its kernel is
  generated by a regular sequence, since $R$ is complete intersection. Proposition
  \eqref{vsmall:regular} implies that $\wh R\otimes_RM$ is \vsmall\, over $Q$, and then
  Theorem \eqref{vsmall:cimap} implies that $\wh R\otimes_RM$ is \vsmall\, over $\wh R$.
  Applying Corollary \eqref{vsmall:complete} one obtains that $\wh R\otimes_RM$ is
  \vsmall\, over $R$. It remains to recall that $M\simeq \wh R\otimes_RM$.
 \end{proof}

\begin{remark}
  The class of \vsmall\, complexes identified by the preceding result is by no means
  exhaustive.  For instance, one can construct many more of those via Proposition
  \eqref{vsmall:completions}, or by mimicking Example \eqref{psmall:sums}.
\end{remark}

Theorem above and Corollary \eqref{psmall:fdim} suggest the following considerations.

\begin{remarks}
  A homologically finite complex $M$ of $R$-modules is said have \emph{finite
    CI-dimension} if there is a diagram of local homomorphisms $R\to R'\gets Q$ such that
\begin{enumerate}[\quad\rm(i)]
\item $R\to R'$ is flat;
\item $Q\to R'$ is surjective, with kernel generated by a regular sequence, and
\item $\fdim_Q(R'\otimes_RM)$ is finite.
\end{enumerate}
Small complexes have this property, as do homologically finite complexes over complete
intersections; see Avramov, Gasharov, and Peeva~\cite{AGP} for details. The question
naturally arises: \emph{is a complex of finite CI-dimension \vsmall?}

Assume $M$ is a complex of finite CI-dimension. Since the complex $R'\otimes_RM$ is
homologically finite over $R'$, and hence also over $Q$, condition (iii) implies that it
is \perf\, over $Q$. Thus, condition (ii) and Theorem \eqref{vsmall:cimap} yield that
$R'\otimes_RM$ is \vsmall\, over $R'$. Thus, what remains before one can conclude that $M$
is \vsmall, is to deal with the

\begin{Question}
  Let $M$ and $N$ be homologically finite complexes of $R$-modules.
\[
\text{Does} \quad (M\dtensor RR')\rbuilds{R'} (N\dtensor RR') \quad \text{imply}\quad
M\rbuilds R N\,?
\]
In words: if $N\dtensor RR'$ is finitely built from $M\dtensor RR'$, then is $N$ finitely
built from $M$?
\end{Question}
\end{remarks}

Theorem \eqref{vsmall:ciring} allows one to specialize many results in Section
\ref{section:testobjects} to the case where $R$ is complete intersection; these throw new
light on homological algebra over such rings. We spell out the one corresponding to
Theorem \eqref{test:asymptotes} because it identifies new scenarios when the complete
intersection property descends along a local homomorphism.

\begin{theorem}
\label{ci:descent}
Let $\psi\col (Q,\fn,l)\to (R,\fm,k)$ be a local homomorphism and $M$ a homologically
finite complex of $R$-modules with $\hh M\ne0$. If $R$ is complete intersection, then
\[
\cxy Qh\leq \codim R + \cxy {\psi} M \qnd \curv Qh\leq \max\{1,\curv {\psi}M\}
\]
In particular, if $\curv \psi M\leq 1$, then $Q$ is complete intersection.  If $\fdim_QM$
is finite, then $Q$ is complete intersection and $\codim Q\leq \codim R$.
\end{theorem}

\begin{remark}
  The theorem above is a common generalization of Parts (1) and (2) of \cite[(13.1)]{AIM},
  which in turn contain the descent of the complete intersection property along flat
  homomorphisms.
\end{remark}

For use in the proof of the theorem, we recall a characterization of complete intersection
rings in terms of asymptotic invariants of $k$.

\begin{chunk}
\label{ci:characterizations}
Let $(R,\fm,k)$ be a local ring. The following conditions are equivalent:
\begin{enumerate}[\quad\rm(i)]
\item $R$ is complete intersection;
\item $\cxy Rk=\codim R$;
\item $\curv Rk\leq 1$.
\end{enumerate}

The proof that (i) $\implies$ (ii) is contained in a result of Tate, while the reverse
implication is due to Gulliksen. The remaining non-trivial implication, (iii) $\implies$
(ii) was established by Avramov. The reader may refer to \cite[\S(8.1)]{Av:barca} for
details.
\end{chunk}

\begin{proof}[Proof of Theorem \emph{\eqref{ci:descent}}]
  From Theorem \eqref{vsmall:ciring} one obtains that $M$ is \vsmall; with this on hand,
  Theorem \eqref{test:asymptotes} with $Y$ specialized to $k$ yields
\[
\cxy \psi k\leq \cxy Rk + \cxy {\psi} M \qnd \curv \psi k\leq \max\{\curv Rk,\curv
{\psi}M\}
\]
It can be deduced easily from \cite[(5.1.2)]{AIM} that $\cxy \psi k=\cxy Qh$ and
$\curv\psi k=\curv Qh$; these, in conjunction with the inequalities above, and
characterizations of complete intersection in \eqref{ci:characterizations}, yield the
inequalities we seek.

If $\curv\psi M\leq 1$, then the inequality $\curv Qh\leq \max\{1,\curv {\psi}M\}$, which
has been verified, implies $\curv Qh\leq 1$, so $Q$ is complete intersection, by
\eqref{ci:characterizations}.

Suppose $\fdim_QM$ is finite. Then $\cxy\psi M=0=\curv\psi M$, by \cite[(7.1.3.1)]{AIM},
so the inequalities that have been established yield $\cxy Qh\leq \cxy Rk$ and $\curv
Qh\leq 1$.  Now we again invoke \eqref{ci:characterizations} to deduce from the latter
inequality that $Q$ is complete intersection, and hence from the former that $\codim Q\leq
\codim R$.
\end{proof}

In the remainder of this article we focus on the problem below, which seeks a converse to
Theorem \eqref{vsmall:ciring}. A positive answer would provide, in conjunction with
Theorem \eqref{vsmall:ciring}, a new homotopical characterization of local complete intersections.

\begin{question}
\label{cichar:question}
Over a local ring $R$, if each homologically finite complex is \vsmall, then is $R$
complete intersection?
\end{question}

The next result takes us a part of the way toward an affirmative answer.

\begin{theorem}
\label{cichar:rings}
Let $R$ be a local ring. If each homologically finite complex of $R$-modules is \vsmall,
then $R$ is Gorenstein.
\end{theorem}

\begin{proof}
  Theorem \eqref{gor:char} yields the desired conclusion, once we observe that $R$ has a
  homologically finite complex of finite injective dimension.  To verify the latter claim,
  let $K$ be the Koszul complex of $R$ and $E$ the injective hull of the residue field $k$
  of $R$. The complex $M=\Hom RKE$ has the desired properties. Indeed, $\hh M$ is non-zero and
  artinian, for example, by Matlis duality, and the maximal ideal of $R$ annihilates
  it. Thus $\hh M$ is a non-zero $k$-vector space of finite rank, and hence a finite
  $R$-module. Since $E$ is injective and $K$ is a finite free complex, $\id_RM$ is finite,
  as desired.
\end{proof}

One may also deduce the preceding result from a version for homomorphisms, described
below; however, the direct proof is easier. We note that a \emph{Gorenstein} homomorphism
is one which is quasi-Gorenstein and has finite flat dimension.

\begin{theorem}
  Let $\psi\col Q\to R$ be a module finite local homomorphism with $\fdim_QR$ finite. If
  each homologically finite complex of $R$-modules that is \perf\, over $Q$ is \vsmall\,
  over $R$, then the homomorphism $\psi$ is Gorenstein.
\end{theorem}

\begin{proof}
  Set $D=\rhom QRQ$; via $R$, this acquires a structure of a complex of $R$-modules.
  Since $\psi$ is module finite, $D$ is a dualizing complex for $\psi$, by
  \cite[(5.12)]{AF:fgd}.  Viewed as a complex of $Q$-modules, $R$ is \perf, by
  \eqref{perf=small}, and hence so is $D$. The hypothesis now entails $D$ is \vsmall\,
  over $R$. Now invoke Theorem \eqref{qgor}.
\end{proof}

The converse to Theorem \eqref{cichar:rings} does not hold; in that, there are Gorenstein
rings having complexes, and even modules, that are homologically finite but not \vsmall.
One family of examples is described below.

\begin{example}
\label{nonci:idealization}
Let $Q$ be a Cohen-Macaulay local ring with dualizing module $\omega$.  Set $R=Q\ltimes
\omega$, the trivial extension of $Q$ by $\omega$, and view $Q$ as an $R$-module via the
canonical surjection $R\to Q$. The local ring $R$ is Gorenstein; see \cite[(3.3.6)]{BH}.
However, if $Q$ is \vsmall\, as an $R$-module, then $Q$ is Gorenstein.  

Indeed, if the $R$-module $Q$ is \vsmall, then Corollary \eqref{vsmall:retracts} implies
that $\fdim_QR$ is finite; thus $\fdim_Q\omega$ is finite, so $Q$ is Gorenstein;
see~\cite[(3.4.12)]{La:gd}, or use (c) $\implies$ (a) of Theorem \eqref{gor:char}, keeping
in mind that $\idim_Q\omega$ is finite.
\end{example}

Now we describe another class of Gorenstein rings having modules that are not \vsmall.
The special case $\rank_k V = 3$ was the first ring recognized, by Macaulay it seems, to
be Gorenstein but not complete intersection; this is one reason for our interest in this
example.

\begin{example}
\label{nonci:macaulay}
Let $k$ be a field of characteristic $\ne 2$ and $V$ a $k$-vector space of finite rank.
Each symmetric bilinear form $\bsq$ on $V$ gives rise to a commutative, graded $k$-algebra
\[
R(\bsq) = k\oplus V \oplus k \with \text{$v\cdot w=q(v,w)$ for $v,w$ in $V$.}
\]
This $k$-algebra is concentrated in degrees $0$, $1$ and $2$, is artinian and local.

\begin{Claim}
  Assume that the form $\bsq$ is \emph{non-degenerate}. The following statements hold:
\begin{enumerate}[\quad\rm(a)]
\item $R(\bsq)$ is Gorenstein; it is complete intersection if and only if $\rank_kV\leq
  2$.
\item When $\rank_kV\geq 3$, the ring $R(\bsq)$ has a module that is not \vsmall.
\end{enumerate}
\end{Claim}

Property (a) is well-known, and one justification for it is given below.  To start with,
it is convenient to set $d=\rank_k V$, and abbreviate $R(\bsq)$ to $R$.  The claims are
all easy to verify when $d=1$, for then $R\cong k[x]/(x^3)$.  Assume henceforth that
$d\geq 2$; choosing an orthogonal basis for $V$ allows one to realize $R$ as
\[
R = \frac{k[x_1,\cdots,x_d]}{J} \where J=\bigl(\{x_1^2-c_ix_i^2\}_{2\les i\les d}\,,
\{x_ix_j\}_{1\les i < j\leq d}\bigr)
\]
where the $c_i$ are non-zero elements in $k$. Note that the prescribed generating set for
$J$ is minimal, and that it contains $(d^2-d+2)/2$ elements.

\smallskip

\emph{Proof of\, \rm{(a)}}: As $\bsq$ is non-degenerate, the socle of $R$ is precisely its
degree two component, and hence of rank one. Since $R$ is artinian, one obtains that it is
Gorenstein, by \cite[(3.2.10)]{BH}. Also, $R$ is complete intersection precisely when the
number of generators of $J$, that is to say, $(d^2-d+2)/2$, equals $d$. This happens if
and only if $d\leq 2$.

\smallskip

\emph{Proof of\, \rm{(b)}}: We claim that the $R$-module below is not \vsmall:
\[
M = R/(x_1x_2,x_1+c_2x_2,x_3,\dots,x_d)
\]
To see this, set $Q= k[x_1,\dots,x_d]/(x_1x_2)$ and let $\psi\col Q\to R$ be the canonical
surjection. By Theorem \eqref{test:fdim}, it suffices to prove that $\fdim_Q M$ is finite,
while $\fdim_QR$ is not. The first claim is easily settled: the sequence
$x_1+c_2x_2,x_3,\dots,x_d$ is regular on $Q$, so the Koszul complex on it is a minimal
resolution of $M$ over $Q$. The non-finiteness of $\fdim_QR$ is the special case
$I=(x_1x_2)$ of the more general statement:

\begin{Claim}
  Let $I\subset J$ be a homogeneous ideal of $k[x_1,\dots,x_d]$; set
  $Q=k[x_1,\dots,x_d]/I$.  If $I$ is generated by a regular sequence, then
  $\fdim_QR=\infty$.
\end{Claim}

Indeed, suppose $n=\fdim_QR$ is finite; then $n=\dim Q$, but this is not important for
what follows. The homogeneous minimal $Q$-free resolution of $R$ reads
\[
0\to \bigoplus_{j\in\BZ}Q(-j)^{b_{nj}}\to \cdots\to \bigoplus_{j\in\BZ}Q(-j)^{b_{1j}}\to
Q\to 0
\]
where the $b_{ij}$ are non-negative integers, and for each $i$, one has $b_{ij}=0$ for all
but finitely many $j$.  From the resolution above, and the additivity of Hilbert series on
short exact sequences, one obtains an equality of formal power series
\begin{equation*}
  \hilb Rs = \hilb Qs \cdot \bigl(\,\sum_{i,j} (-1)^i b_{ij}s^j\bigr)\tag{$\dagger$}
\end{equation*}
Let $r_1,\dots,r_c$ be a regular sequence generating $I$ and set $d_i=\deg{r_i}$;
note that $d_i\geq 2$, since $I\subset J$.  
\[
\hilb Rs = 1+ds+s^2 \qnd \hilb Qs = \frac{\prod_{i=1}^c(1-t^{d_i})}{(1-t)^d}
\]
These formulas combined with ($\dagger$) imply that the primitive $d_1$th root of unity is
a root of $1+ds+s^2$. This cannot happen because $d\geq 3$ and $d_1\geq 2$, 
\end{example}

\section*{Acknowledgments}
It is a pleasure to thank Lucho Avramov for numerous discussions on this work, and Sean
Sather-Wagstaff for comments on this manuscript.  Part of this research was done at the
Mathematical Sciences Research Institute, Berkeley, where the authors met in February 2003
during the special year in commutative algebra.  We are grateful to the MSRI for its
hospitality. The second author thanks the University of Nebraska, Lincoln, for hospitality
during a visit in November 2003.

\end{document}